\documentclass[a4paper,11pt,oneside]{amsart}

\pdfoutput=1

\usepackage[english]{babel}

\usepackage{lmodern}
\usepackage[T1]{fontenc}
\usepackage{microtype}
\usepackage{typearea}

\usepackage{amsmath}
\usepackage{amsthm}
\usepackage{amssymb}
\usepackage{mathtools}
\usepackage{mathrsfs}
\usepackage{dsfont}

\usepackage{calc}
\usepackage{stmaryrd}
\usepackage{stackrel}
\usepackage{colonequals}
\usepackage{chngcntr}
\usepackage{multirow}
\usepackage{enumerate} 
\usepackage[autostyle]{csquotes}
\usepackage[all]{xy}
\usepackage{tikz,pgffor}
\usepackage[outline]{contour}
\contourlength{2pt}

\allowdisplaybreaks
\numberwithin{equation}{section}

\newcommand{\rleft}{\mathopen{}\mathclose\bgroup\left}
\newcommand{\rright}{\aftergroup\egroup\right}

\newtheorem{theorem}{Theorem}[section]
\newtheorem{lemma}[theorem]{Lemma}
\newtheorem{prop}[theorem]{Proposition}

\newtheorem{cor}[theorem]{Corollary}
\theoremstyle{definition}
\newtheorem{definition}[theorem]{Definition}
\newtheorem{remark}[theorem]{Remark}
\newtheorem{example}[theorem]{Example}

\numberwithin{equation}{section}

\DeclareMathOperator*{\Cl}{Cl}

\DeclareMathOperator*{\cone}{cone}

\DeclareMathOperator*{\Spec}{Spec}
\DeclareMathOperator*{\Hom}{Hom}

\DeclareMathOperator{\vspan}{span}

\newcommand{\SL}{\mathrm{SL}}

\newcommand{\ie}{i.\,e.~}

\newcommand{\C}{\mathds{K}}
\newcommand{\Q}{\mathds{Q}}

\newcommand{\V}{\mathds{V}}

\newcommand{\Z}{\mathds{Z}}

\newcommand{\Om}{\mathcal{O}}
\newcommand{\Qm}{\mathcal{Q}}
\newcommand{\Vm}{\mathcal{V}}
\newcommand{\Cm}{\mathcal{C}}
\newcommand{\Fm}{\mathcal{F}}
\newcommand{\Dm}{\mathcal{D}}

\newcommand{\Rm}{\mathcal{R}}
\newcommand{\Pm}{\mathcal{P}}

\newcommand{\Nm}{\mathcal{N}}
\newcommand{\Mm}{\mathcal{M}}

\newcommand{\Df}{\mathcal{D}}

\newcommand{\Ff}{\mathfrak{F}}
\newcommand{\Af}{\mathfrak{A}}
\newcommand{\Bf}{\mathfrak{B}}

\newcommand{\adiv}{\Delta}
\newcommand{\bdiv}{\Dm}
\newcommand{\gdiv}{\Gamma}

\begin{document}

\title{Spherical varieties with the $A_k$-property}

\author{Giuliano Gagliardi}
\address{Institut f\"ur Algebra, Zahlentheorie und Diskrete Mathematik, Leibniz Universit\"at Hannover,
Welfengarten 1, 30167 Hannover, Germany}
\email{gagliardi@math.uni-hannover.de}


\begin{abstract}
An algebraic variety is said to have the $A_k$-property
if any $k$ points are contained in some common affine open neighbourhood.
A theorem of W{\l}odarczyk states that a normal variety
has the $A_2$-property
if and only if it admits a closed
embedding into a toric variety.
Spherical varieties can be regarded as a generalization of toric
varieties, but they do not have the $A_2$-property in general.
We provide a combinatorial criterion for the
$A_k$-property of spherical varieties 
by combining the theory of bunched rings with the
Luna-Vust theory of spherical embeddings.
\end{abstract}

\maketitle

\section{Introduction}
\label{sec:int}
Throughout the paper,
we work with algebraic varieties and algebraic groups
over an algebraically closed field $\C$ of characteristic zero.

\begin{definition}[{see, for instance, \cite{a2,swi01,hau02}}]
A variety $X$ is said to have the \emph{$A_k$-property}
if any $k$ points $x_1, \dots, x_k \in X$ are contained in some common affine open
neighbourhood.
\end{definition}
Clearly, any quasi-projective variety has the $A_k$-property
for every $k$. According to the generalized 
Kleiman-Chevalley criterion for quasi-projectivity (see~\cite{kle66, wlo99, ben13}),
the converse is true for normal varieties.

There exist toric varieties of dimension $3$ and greater
which are not quasi-projective
(see~\cite[after~2.16]{oda}),
but they always have the $A_2$-property.
In fact, W{\l}odarczyk has shown in \cite{a2} that a normal variety has the $A_2$-property
if and only if it admits a closed embedding
into a toric variety.

In this paper, we consider spherical varieties, which may
be regarded as a generalization
of toric varieties. Fix a connected reductive group $G$
and a Borel subgroup $B \subseteq G$.
A closed subgroup $H \subseteq G$ is
called \emph{spherical} if $G/H$ contains an open $B$-orbit,
and then $G/H$ is called a \emph{spherical homogeneous space}.
A $G$-equivariant open embedding $G/H \hookrightarrow X$ into
a normal irreducible $G$-variety $X$ is called a \emph{spherical embedding}, 
and then $X$ is called a \emph{spherical variety}.

According to the Luna-Vust theory (see \cite{lunavust, knopsph}), 
we can associate to any
spherical embedding $G/H \hookrightarrow X$ a combinatorial
object called a \emph{colored fan}.
We denote by $\Mm$ the weight lattice of $B$-semi-invariants
in the function field $\C(G/H)$ and by $\Nm_\Q$
the vector space dual to $\Mm_\Q \coloneqq \Mm \otimes_\Z \Q$.
We denote by $\adiv$ the set
of $B$-invariant prime divisors in $X$.
The subset $\bdiv \subseteq \adiv$ of $B$-invariant prime divisors in $G/H$ is
called the set of \emph{colors}.
Moreover, there is a natural map $\rho\colon \adiv \to \Nm_\Q$.
The colored fan associated to the spherical embedding $G/H \hookrightarrow X$
is then given by
\begin{align*}
\Sigma \coloneqq \rleft\{(\cone(\rho(I_Y)), I_Y \cap \Dm) : \text{$Y \subseteq X$ is a $G$-orbit, } I_Y = \{D \in \adiv : Y \subseteq \overline{D}\}\rright\}\text{,}
\end{align*}
which means that the colored fan $\Sigma$ contains a pair $(\Cm, \Fm)$, called a \emph{colored cone},
for every $G$-orbit $Y \subseteq X$.
For details, we refer to Section~\ref{sec:cf}.

It was noticed by Huruguen that, in contrast to toric varieties, spherical
varieties do not have the $A_2$-property in general.
In fact, according to \cite[Remark~2.38]{hur11},
the example of a non-projective spherical variety with $\dim \Nm_\Q = 2$
considered in \cite[Remarque~3.11]{pauer} and \cite[Example~17.7]{ti}
fails to have the $A_2$-property.
It is a spherical embedding of $\SL_3(\C)/\SL_2(\C)$,
whose colored fan is shown in the following picture. This example can
be generalized to higher dimensions (see~\cite[Example~4.2]{g1}).

\begin{align*}
\begin{tikzpicture}[scale=0.8]
	\clip (-1.5, -1.5) -- (1.5, -1.5) -- (1.5, 1.5) -- (-1.5, 1.5) -- cycle;
	\fill[color=gray!30] (-3, 3) -- (3, -3) -- (-3, -3) -- cycle;
	\draw (1, 0) circle (2pt);
	\draw (0, 1) circle (2pt);
	\draw (0,0) -- (3,0);
	\draw (0,0) -- (0,3);
	\draw[very thick] (0,0) -- (-4, 2);
	\draw[very thick] (0,0) -- (2, -4);
        \fill (0,9pt) circle (2pt);
        \fill (13pt,0) circle (2pt);
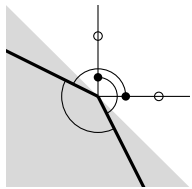
	\begin{scope}\clip (0,0) -- (1, -2) -- (2, 0) -- (0, 2) -- cycle; \draw (0,0) circle (9pt);\end{scope}
	\begin{scope}\clip (0,0) -- (2, 0) -- (0, 2) -- (-2, 1) -- cycle; \draw (0,0) circle (13pt);\end{scope}
	\begin{scope}\clip (0,0) -- (1, -2) -- (-2, -2) -- (-2, 1) -- cycle; \draw (0,0) circle (17pt);\end{scope}
\end{tikzpicture}
\end{align*}
Note that two of the colored cones in this colored fan do not intersect in a common
face (we refer to Section~\ref{sec:cf} for the precise definition of \enquote{face}),
which is allowed by the Luna-Vust theory
as long as this does not happen inside a certain valuation cone $\Vm \subseteq \Nm_\Q$, which is shown in grey.

In order to give a precise characterization of
the $A_k$-property, we define
\begin{align*}
\Sigma^\sharp \coloneqq \{\cone([\adiv \setminus I]) : I \subseteq \adiv, (\cone(\rho(I)), I \cap \bdiv) \in \Sigma\}\text{,}
\end{align*}
which is a set of cones inside
the vector space $\Cl(X)_\Q \coloneqq \Cl(X) \otimes_\Z \Q$.
The divisor class group $\Cl(X)$ is generated by the divisor classes
$[D]$ for $D \in \adiv$, and the relations are given in \cite[Proposition~4.1.1]{brcox}. For details, we refer to Sections~\ref{sec:gdse} and \ref{sec:br}.

If $X$ has the $A_2$-property and $\Gamma(X, \Om_X^*) = \C^*$,
then we will see in Section~\ref{sec:br} that $\Sigma^\sharp$ is the bunch of cones associated to $X$
by the theory of bunched rings, but our main
result also holds without these assumptions.

\begin{theorem}
\label{th:ak}
Let $G/H \hookrightarrow X$ be a spherical
embedding with associated colored fan $\Sigma$.
Then $X$ has the $A_k$-property if and
only if
for any $k$ cones $\tau_1, \dots, \tau_k \in \Sigma^\sharp$ we have $\tau_1^\circ \cap \dots \cap \tau_k^\circ \ne \emptyset$
(where $\tau_i^\circ$ denotes the relative interior of $\tau_i$).
\end{theorem}

For the $A_2$-property,
we obtain the following characterization, which can be verified
on the colored fan $\Sigma$ itself.

\begin{definition}
\label{def:ccint}
The \emph{intersection} of two colored cones $(\Cm_1, \Fm_1)$
and $(\Cm_2, \Fm_2)$ is defined to be the colored cone $(\Cm_1 \cap \Cm_2, \Fm_1 \cap \Fm_2)$.
\end{definition}

\begin{theorem}
\label{th:m}
Let $G/H \hookrightarrow X$ be a spherical
embedding with associated colored fan $\Sigma$.
Then $X$ has the $A_2$-property if and
only if any two colored cones in $\Sigma$
intersect in a common face.
\end{theorem}

A spherical variety is called \emph{horospherical} if $\Vm = \Nm_\Q$.

\begin{cor}
Every horospherical variety has the $A_2$-property.
\end{cor}

\section{Spherical embeddings and colored fans}
\label{sec:cf}

We give a brief overview over the parts of the Luna-Vust theory of spherical embeddings which are relevant for us.
For details, we refer to \cite{lunavust, knopsph}. A survey can also be found in \cite{ti}.

Let $G/H$ be a spherical homogeneous space.
We denote by $\Mm$ the weight lattice of $B$-semi-invariants
in the function field $\C(G/H)$ and by $\Nm \coloneqq \Hom(\Mm, \Z)$ the
dual lattice, together with the
natural pairing \[\langle \cdot, \cdot \rangle\colon \Nm \times \Mm \to \Z\text{.}\]

We denote by $\Dm$
the set of $B$-invariant prime divisors in $G/H$.
The elements in $\Dm$ are called the \emph{colors}.
Moreover, we denote by $\rho: \Dm \to \Nm$ the map given by
$\langle \rho(D), \chi \rangle \coloneqq \nu_D(f_\chi)$ for $D \in \Dm$ where
$\nu_D$ is the discrete valuation on $\C(G/H)$ induced by the prime
divisor $D$ and
$f_\chi \in \C(G/H)$ is a $B$-semi-invariant rational function of weight $\chi \in \Mm$ (such
a rational function $f_\chi$ is uniquely determined up to a constant factor).

In the same way, we define a map $\Vm \to \Nm_\Q$ from the set $\Vm$ of $G$-invariant discrete
valuations on $\C(G/H)$. This map is injective, so that we may
consider $\Vm$ as a subset of $\Nm_\Q$. It is known from \cite{brg} that
$\Vm$ is a cosimplicial (in particular full-dimensional) cone, called the \emph{valuation cone}
of $G/H$.
The objects \[\rho\colon \Dm \to \Nm_\Q \supseteq \Vm\] are called
a \emph{colored vector space}.

\begin{definition}
A \emph{colored cone} is a pair $(\Cm, \Fm)$
such that $\Fm \subseteq \Dm$ is a subset and $\Cm \subseteq \Nm_\Q$ is a cone
generated by $\rho(\Fm)$ and finitely many elements of $\Vm$. 
It is called
\begin{enumerate}[(i)]
\item \emph{supported} if $\Cm^\circ \cap \Vm \ne \emptyset$, where $\Cm^\circ$
denotes the relative interior of $\Cm$,
\item \emph{pointed} if $\Cm$ is pointed
and $0 \notin \rho(\Fm)$,
\item \emph{simplicial} if $\Cm$ is spanned
by a part of a $\Q$-basis of $\Nm_\Q$ which contains $\rho(\Fm)$ and $\rho|_{\Fm}\colon \Fm \to \Nm_\Q$
is injective.
\end{enumerate}
A \emph{face} of a colored cone $(\Cm, \Fm)$ is a colored cone $(\Cm_0, \Fm_0)$
such that $\Cm_0$ is a face of $\Cm$ and $\Fm_0 = \Fm \cap \rho^{-1}(\Cm_0)$.
A \emph{colored fan} is a nonempty set $\Sigma$ of pointed supported
colored cones such that every supported face of a colored cone in $\Sigma$
also belongs to $\Sigma$ and
for every $u \in \Vm$ there is at most one $(\Cm, \Fm) \in \Sigma$ with
$u \in \Cm^\circ$.
\end{definition}

\begin{remark}
In contrast to some of the literature, we do not require colored cones
and their faces to be supported.
\end{remark}

For a spherical embedding $G/H \hookrightarrow X$, we denote
by $\gdiv$ the set of $G$-invariant prime divisors in $X$.
Then $\adiv \coloneqq \bdiv \cup \gdiv$ is the set of all $B$-invariant
prime divisors, and the definition of $\rho\colon \bdiv \to \Nm_\Q$
extends to $\rho\colon \adiv \to \Nm_\Q$.

For any $G$-orbit $Y \subseteq X$, we denote by $I_Y \subseteq \adiv$
the set of $B$-invariant prime divisors containing $Y$ in their closure
and we set \[(\Cm_Y, \Fm_Y) \coloneqq (\cone(\rho(I_Y)), I_Y \cap \Dm)\text{.}\]

\begin{theorem}[{\cite[Theorem~3.3]{knopsph}}]
\label{th:bijcfan}
The map
\begin{align*}
(G/H \hookrightarrow X) \mapsto \Sigma \coloneqq \{(\Cm_Y, \Fm_Y): Y \subseteq X \text{ is a $G$-orbit}\}
\end{align*}
defines a bijection between isomorphism classes of spherical embeddings of $G/H$ and
colored fans.

Moreover, the assignment
\begin{align*}
\{\text{$G$-orbits in $X$}\} &\to \Sigma\\
Y &\mapsto (\Cm_Y, \Fm_Y)
\end{align*}
is a bijection such that for two $G$-orbits $Y_1, Y_2 \subseteq X$
we have $Y_1 \subseteq \overline{Y_2}$ if and only if $(\Cm_{Y_2}, \Fm_{Y_2})$ is a face of $(\Cm_{Y_1}, \Fm_{Y_1})$.
\end{theorem}

\begin{remark}
Note that $\rho|_\bdiv$ need not be injective (see, for instance, Example~\ref{ex:2}).
On the other hand, Theorem~\ref{th:bijcfan} implies
that the elements of $\rho(\gdiv)$ generate pairwise different nonzero
rays. It also implies $\rho(\gdiv) \subseteq \Vm$. Moreover, for 
$D' \in \bdiv$ with $\rho(D') \in \Vm$ it is possible that there
exists $D'' \in \gdiv$ such that $\rho(D')$ and $\rho(D'')$ generate
the same ray.
\end{remark}

It follows from \cite[Theorem~6.6]{knopsph} that,
under the orbit-cone correspondence of Theorem~\ref{th:bijcfan},
the $G$-orbits in $X$ of codimension $1$ correspond to
the colored cones in $\Sigma$ of the form $(\cone(u), \emptyset)$ for $u \in \Nm_\Q$ (which
implies $u \in \Vm$). These are exactly the (pairwise distinct)
colored cones $(\cone(\rho(D)), \emptyset)$ for $D \in \gdiv$.

\begin{remark}[{\cite[Proposition~3.1]{bri89}}]
Let $G/H \hookrightarrow X$ be a spherical embedding with associated colored fan $\Sigma$.
Then $X$ is $\Q$-factorial if and only if every colored cone in $\Sigma$ is simplicial.
\end{remark}

\begin{remark}
\label{rem:intf}
Two colored cones $(\Cm_1, \Fm_1)$ and $(\Cm_2, \Fm_2)$
intersect in a common face if and only if there
exists $e \in \Mm_{\Q}$ with $e|_{\Cm_1} \ge 0$, $e|_{\Cm_2} \le 0$,
\begin{align*}
\Cm_1 \cap e^\perp = \Cm_1 \cap \Cm_2 = e^\perp \cap \Cm_2\text{,} &&\text{and}&&
\Fm_1 \cap \rho|_{\Dm}^{-1}(e^\perp) = \rho|_{\Dm}^{-1}(e^\perp) \cap \Fm_2\text{.}
\end{align*}
\end{remark}

\begin{remark}
\label{rem:constunits}
According to \cite[Proposition~1.3(ii)]{kkv89}, the units in
$\Gamma(X, \Om_X)$ are $G$-semi-invariant.
In particular,
they are $B$-semi-invariant, so that
we  have $\Gamma(X, \Om^*_X) = \C^*$ if and only if the set $\rho(\adiv)$
generates $\Nm_\Q$ as a vector space.
\end{remark}

The following result will be used in Section~\ref{sec:br}.

\begin{prop}
\label{prop:qp}
Let $G/H \hookrightarrow X$ be a spherical embedding with at most one
closed $G$-orbit of codimension at least $2$. Then $X$ is quasi-projective.
\end{prop}
\begin{proof}
Let $Y \subseteq X$ be a $G$-orbit of maximal codimension (this is the unique
closed $G$-orbit of codimension at least $2$ if such an orbit exists).
We denote by $X_Y \subseteq X$ the open $G$-stable subvariety obtained
by removing from $X$ all $G$-orbits not containing $Y$
in their closure.
As $Y$ is the unique closed $G$-orbit in $X_Y$,
it follows from \cite[Lemma~8]{sum1} that $X_Y$ is quasi-projective.

According to \cite[Theorem~4.9]{sum2},
there exists a $G$-equivariant surjective morpishm $q\colon X' \to X$
where $X'$ is a quasi-projective $G$-variety and $q$ is an isomorphism over $X_Y$.
We may moreover assume $X'$ to be normal.
Let $\Sigma$ be the colored fan associated to $G/H \hookrightarrow X$,
and let $\Sigma'$ be the colored fan associated to $G/H \hookrightarrow X'$.
Since the $G$-orbits in $X \setminus X_Y$ are all of codimension $1$, they correspond
to the colored cones $(\cone(\rho(D)), \emptyset) \in \Sigma$ with $D \in \gdiv_{X \setminus X_Y}$ for
some subset $\gdiv_{X\setminus X_Y} \subseteq \gdiv$. Then
\cite[Theorem~4.1]{knopsph} implies that we have
$(\cone(\rho(D)), \emptyset) \in \Sigma'$ for every $D \in \gdiv_{X\setminus X_Y}$,
hence we obtain $X \subseteq X'$.
\end{proof}

\section{Gale duality for spherical embeddings}
\label{sec:gdse}

The aim of this section is to generalize 
some results on Gale duality for toric varieties,
as presented in \cite[2.2.1]{coxrings},
to the setting of spherical embeddings.

There are two main differences
between the toric and the spherical case. The first difference is that
there could exist distinct $D, D' \in \adiv$ such that
$\rho(D)$ and $\rho(D')$ generate the same ray in $\Nm_\Q$
when colors are involved.
The reason why this is not a problem is that colored cones by definition
keep track of the colors.
The second difference is that in the case $\Vm \ne \Nm_\Q$ some cones
will be ignored when they are not supported according to the definitions below.

Let $\rho\colon \bdiv \to \Nm_\Q \supseteq \Vm$ be a colored vector space,
and let $\gdiv$ be a finite set equipped with another map
$\rho\colon \Gamma \to \Vm$.
We define $\Delta \coloneqq \Dm \cup \Gamma$ to be a disjoint union
and obtain a map $\rho\colon \adiv \to \Nm_{\Q}$.
We assume $\vspan_\Q{} \rho(\Delta) = \Nm_\Q$.
Let $\Q^{\adiv}$ and $(\Q^{\adiv})^*$ be dual vector
spaces with respective standard bases $\{e_D : D \in \adiv\}$ and
$\{e_D^* : D \in \adiv\}$, which are dual to each other.
The map $\Pm\colon (\Q^{\adiv})^* \to \Nm_\Q$ with $e_D^* \mapsto \rho(D)$
induces the following pair of mutually dual exact sequences of vector spaces.
\begin{align*}
\xymatrix@R=0pt{
0 \ar[r] & L_\Q \ar[r] & {(\Q^{\adiv})^*} \ar[r]^-{\Pm} & \Nm_\Q \ar[r] & 0\\
0 & K_\Q \ar[l] & {\phantom{(}\Q^{\adiv}\phantom{)^*}} \ar[l]_-{\Qm} & \Mm_\Q \ar[l] & 0 \ar[l]
}
\end{align*}
For $D \in \adiv$, we also write $\Pm(D)$ for $\Pm(e^*_D)$ and $\Qm(D)$ for $\Qm(e_D)$.

\begin{definition}
A \emph{$\Pm$-cone} is a pair $(\cone(\Pm(I)), I \cap \bdiv)$ where
$I \subseteq \adiv$.
\end{definition}

\begin{remark}
Note that every $\Pm$-cone is a colored cone.
In particular, the definitions of \enquote{supported},
intersections, faces, \enquote{pointed}, and
\enquote{simplicial} from Section~\ref{sec:cf}
are applicable.
\end{remark}

\begin{definition}
A \emph{$\Qm$-cone} is a cone in $K_\Q$ generated by a subset of $\Qm(\adiv)$.
\end{definition}

\begin{definition}
For any set $\Sigma$ of $\Pm$-cones and any set
$\Theta$ of $\Qm$-cones we define
\begin{align*}
\Sigma^\natural &\coloneqq \{\cone(\Qm(\adiv \setminus I)) : I \subseteq \adiv,
(\cone(\Pm(I)), I \cap \bdiv) \in \Sigma\}\text{,}\\
\Theta^\natural &\coloneqq \{(\cone(\Pm(\adiv \setminus J)), (\adiv \setminus J) \cap \bdiv) : J \subseteq \adiv, \cone(\Qm(J)) \in \Theta\}\text{.}
\end{align*}
\end{definition}

\begin{definition}
A $\Qm$-cone $\tau$ is called \emph{supported} if $\{\tau\}^\natural$ contains a supported $\Pm$-cone.
For any set $\Sigma$ of $\Pm$-cones we define
\begin{align*}
\overline{\Sigma} \coloneqq \{(\Cm, \Fm) \in \Sigma : \text{$(\Cm, \Fm)$ is supported}\}\text{,}
\end{align*}
\ie the $\Pm$-cones which are not supported are removed from $\Sigma$.
\end{definition}

\begin{remark}
A $\Qm$-cone $\tau$ is supported if and only if $\overline{\{\tau\}^\natural}$ is not empty.
\end{remark}

\begin{definition}
A set $\Sigma$ of $\Pm$-cones is called a \emph{$\Pm$-quasifan} if it is nonempty and
\begin{enumerate}[(i)]
\item every $(\Cm, \Fm) \in \Sigma$ is supported,
\item any two $(\Cm_1, \Fm_1), (\Cm_2, \Fm_2) \in \Sigma$ intersect in a common face,
\item for any $(\Cm, \Fm) \in \Sigma$ every supported face of $(\Cm, \Fm)$ also belongs to $\Sigma$.
\end{enumerate}
A \emph{$\Pm$-fan} is a $\Pm$-quasifan consisting
of pointed $\Pm$-cones.
A $\Pm$-(quasi)fan is called \emph{maximal} if it cannot be extended by adding supported $\Pm$-cones.
It is called \emph{true} if it contains the $\Pm$-cone $(0, \emptyset)$ in the case $\Df \ne \emptyset$
as well as the $\Pm$-cones $(\cone(\rho(D)), \emptyset)$ for $D \in \gdiv$.
\end{definition}

\begin{definition}
A set $\Theta$ of $\Qm$-cones is called a \emph{$\Qm$-bunch} if it is nonempty and
\begin{enumerate}[(i)]
\item every $\tau \in \Theta$ is supported,
\item for any $\tau_1, \tau_2 \in \Theta$ we have $\tau_1^\circ \cap \tau_2^\circ \ne \emptyset$,
\item for any $\tau \in \Theta$ every supported $\Qm$-cone $\tau_0$ with $\tau^\circ \subseteq \tau_0^\circ$ 
also belongs to $\Theta$.
\end{enumerate}
\item A $\Qm$-bunch $\Theta$ is called \emph{maximal} if it cannot be extended by adding supported $\Qm$-cones.
It is called \emph{true} if it contains the $\Qm$-cone $\cone(\Qm(\adiv))$ in the case
$\Df \ne \emptyset$ as well as 
the $\Qm$-cones $\cone(\Qm(\adiv \setminus \{D\}))$ for $D \in \gdiv$.
\end{definition}

We can now state our generalization of \cite[Theorem~2.2.1.14]{coxrings} to
the spherical situation. The proof will be given in Section~\ref{sec:p}.

\begin{theorem}\label{thm:fb}
We have an order reversing map
\begin{align*}
\{\text{$\Qm$-bunches}\} \to \{\text{$\Pm$-quasifans}\}, \qquad
\Theta \mapsto \Theta^{\sharp} \coloneqq \overline{\Theta^\natural}\text{.}
\end{align*}
Now assume that the elements in $\rho(\gdiv)$ generate pairwise different rays.
Then there are mutually inverse order reversing bijections
\begin{align*}
\{\text{true maximal $\Qm$-bunches}\} &\leftrightarrow \{\text{true maximal $\Pm$-fans}\}\text{,}\\
\Theta &\mapsto \Theta^{\sharp}\text{,}\\
\Sigma^\natural \eqqcolon \Sigma^{\sharp} &\mapsfrom \Sigma\text{.}
\end{align*}
Under these bijections, the true maximal $\Pm$-fans consisting of simplicial cones correspond to
the true maximal $\Qm$-bunches consisting of full-dimensional cones.
\end{theorem}

\begin{example}
\label{ex:1}
We consider the spherical homogeneous space $\SL_3(\C)/\SL_2(\C)$.
The left-hand picture below shows the vector space $\Nm_\Q$
with the valuation cone $\Vm$ in grey. There
are two colors, $\Dm = \{D_1, D_2\}$, where the elements
$\rho(D_i) \in \Nm_\Q$ are represented by a circle.
For details, we refer to \cite{pauer}, \cite[Example~17.7]{ti},
or \cite[Example~4.2]{g1}.
We have added two $\SL_3(\C)$-invariant prime divisors,
$\gdiv = \{D_3, D_4\}$. In the picture, the number $i$ is written near the ray generated by $\Pm(D_i)$.
The right-hand picture shows the vector space $K_\Q$ and the number $i$ is written near the ray generated by $\Qm(D_i)$.
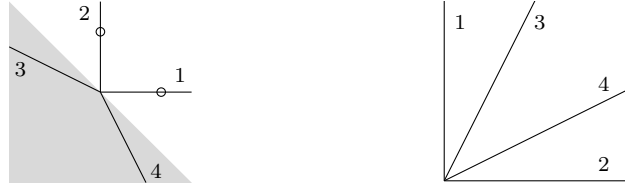
\begin{align*}
\begin{tikzpicture}[scale=0.8]
	\clip (-1.5, -1.5) -- (1.5, -1.5) -- (1.5, 1.5) -- (-1.5, 1.5) -- cycle;
	\fill[color=gray!30] (-3, 3) -- (3, -3) -- (-3, -3) -- cycle;
	\draw (1, 0) circle (2pt);
	\draw (0, 1) circle (2pt);
	\draw (0,0) -- (1.5,0) node[very near end,above]{${\scriptstyle 1}$};
	\draw (0,0) -- (0,1.5) node[very near end,left]{${\scriptstyle 2}$};
	\draw (0,0) -- (-1.5, 0.75) node[very near end,below]{${\scriptstyle 3}$};
	\draw (0,0) -- (0.75, -1.5) node[very near end,right]{${\scriptstyle 4}$};
\end{tikzpicture}
&&\begin{tikzpicture}[scale=0.8]
	\clip (-0.03, -0.03) -- (-0.03, 2.97) -- (2.97, 2.97) -- (2.97, -0.03) -- cycle;
	\draw (0,0) -- (3, 0) node[very near end,above]{${\scriptstyle 2}$};
	\draw (0,0) -- (0, 3) node[very near end,right]{${\scriptstyle 1}$};;
        \draw (0,0) -- (3, 1.5) node[very near end,above]{${\scriptstyle 4}$};
        \draw (0,0) -- (1.5, 3) node[very near end,right]{${\scriptstyle 3}$};
\end{tikzpicture}
\end{align*}
The following pictures show all possible true maximal $\Pm$-fans $\Sigma_j$ and their corresponding true maximal $\Qm$-bunches $\Theta_j$. In the pictures, the included cones of dimension $2$ are represented by arcs 
while the included cones of dimension $1$ are represented by thick rays. For colored cones,
the colors which are included are indicated by large black dots.
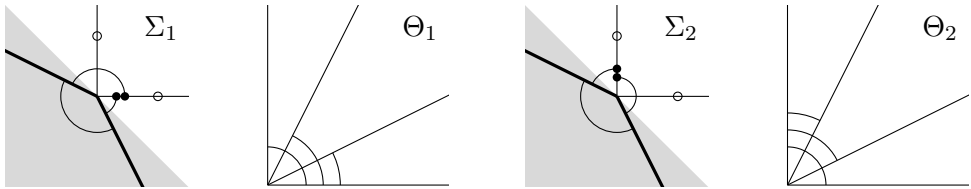
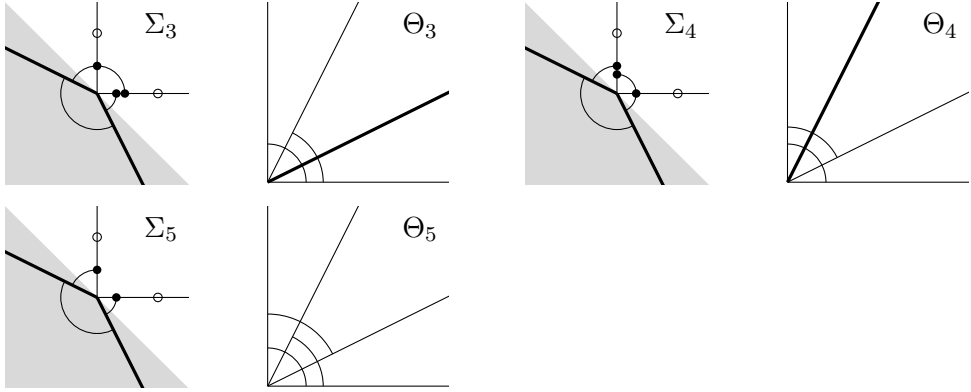
\begin{align*}
\begin{tikzpicture}[scale=0.8]
	\clip (-1.5, -1.5) -- (1.5, -1.5) -- (1.5, 1.5) -- (-1.5, 1.5) -- cycle;
	\fill[color=gray!30] (-3, 3) -- (3, -3) -- (-3, -3) -- cycle;
	\draw (1, 0) circle (2pt);
	\draw (0, 1) circle (2pt);
	\draw (0,0) -- (3,0);
	\draw (0,0) -- (0,3);
	\draw[very thick] (0,0) -- (-4, 2);
	\draw[very thick] (0,0) -- (2, -4);
        \fill (9pt,0) circle (2pt);
        \fill (13pt,0) circle (2pt);
	\begin{scope}\clip (0,0) -- (1, -2) -- (2, 0) -- cycle; \draw (0,0) circle (9pt);\end{scope}
	\begin{scope}\clip (0,0) -- (2, 0) -- (0, 2) -- (-2, 1) -- cycle; \draw (0,0) circle (13pt);\end{scope}
	\begin{scope}\clip (0,0) -- (1, -2) -- (-2, -2) -- (-2, 1) -- cycle; \draw (0,0) circle (17pt);\end{scope}
	\node[anchor=north east] at (1.5, 1.5) {$\Sigma_1$};
\end{tikzpicture}
&&\begin{tikzpicture}[scale=0.8]
	\clip (-0.03, -0.03) -- (-0.03, 2.97) -- (2.97, 2.97) -- (2.97, -0.03) -- cycle;
	\draw (0,0) -- (3, 0);
	\draw (0,0) -- (0, 3);
        \draw (0,0) -- (6, 3);
        \draw (0,0) -- (3, 6);
	\begin{scope}\clip (0,0) -- (3, 0) -- (0, 3) -- cycle; \draw (0,0) circle (18pt);\end{scope}
        \begin{scope}\clip (0,0) -- (3, 0) -- (3, 6) -- cycle; \draw (0,0) circle (26pt);\end{scope}
        \begin{scope}\clip (0,0) -- (3, 0) -- (6, 3) -- cycle; \draw (0,0) circle (34pt);\end{scope}
	\node[anchor=north east] at (2.97, 2.97) {$\Theta_1$};
\end{tikzpicture}
&&
\begin{tikzpicture}[scale=0.8]
	\clip (-1.5, -1.5) -- (1.5, -1.5) -- (1.5, 1.5) -- (-1.5, 1.5) -- cycle;
	\fill[color=gray!30] (-3, 3) -- (3, -3) -- (-3, -3) -- cycle;
	\draw (1, 0) circle (2pt);
	\draw (0, 1) circle (2pt);
	\draw (0,0) -- (3,0);
	\draw (0,0) -- (0,3);
	\draw[very thick] (0,0) -- (-4, 2);
	\draw[very thick] (0,0) -- (2, -4);
        \fill (0,9pt) circle (2pt);
        \fill (0,13pt) circle (2pt);
	\begin{scope}\clip (0,0) -- (1, -2) -- (2, 0) -- (0, 2) -- cycle; \draw (0,0) circle (9pt);\end{scope}
	\begin{scope}\clip (0,0) -- (0, 2) -- (-2, 1) -- cycle; \draw (0,0) circle (13pt);\end{scope}
	\begin{scope}\clip (0,0) -- (1, -2) -- (-2, -2) -- (-2, 1) -- cycle; \draw (0,0) circle (17pt);\end{scope}
	\node[anchor=north east] at (1.5, 1.5) {$\Sigma_2$};
\end{tikzpicture}
&&\begin{tikzpicture}[scale=0.8]
	\clip (-0.03, -0.03) -- (-0.03, 2.97) -- (2.97, 2.97) -- (2.97, -0.03) -- cycle;
	\draw (0,0) -- (3, 0);
	\draw (0,0) -- (0, 3);
        \draw (0,0) -- (6, 3);
        \draw (0,0) -- (3, 6);
	\begin{scope}\clip (0,0) -- (3, 0) -- (0, 3) -- cycle; \draw (0,0) circle (18pt);\end{scope}
        \begin{scope}\clip (0,0) -- (6, 3) -- (0, 6) -- cycle; \draw (0,0) circle (26pt);\end{scope}
        \begin{scope}\clip (0,0) -- (3, 6) -- (0, 3) -- cycle; \draw (0,0) circle (34pt);\end{scope}
	\node[anchor=north east] at (2.97, 2.97) {$\Theta_2$};
\end{tikzpicture}
\\
\begin{tikzpicture}[scale=0.8]
	\clip (-1.5, -1.5) -- (1.5, -1.5) -- (1.5, 1.5) -- (-1.5, 1.5) -- cycle;
	\fill[color=gray!30] (-3, 3) -- (3, -3) -- (-3, -3) -- cycle;
	\draw (1, 0) circle (2pt);
	\draw (0, 1) circle (2pt);
	\draw (0,0) -- (3,0);
	\draw (0,0) -- (0,3);
	\draw[very thick]  (0,0) -- (-4, 2);
	\draw[very thick]  (0,0) -- (2, -4);
        \fill (9pt,0) circle (2pt);
        \fill (13pt,0) circle (2pt);
	\fill (0,13pt) circle (2pt);
	\begin{scope}\clip (0,0) -- (1, -2) -- (2, 0) -- cycle; \draw (0,0) circle (9pt);\end{scope}
	\begin{scope}\clip (0,0) -- (2, 0) -- (0, 2) -- (-2, 1) -- cycle; \draw (0,0) circle (13pt);\end{scope}
	\begin{scope}\clip (0,0) -- (1, -2) -- (-2, -2) -- (-2, 1) -- cycle; \draw (0,0) circle (17pt);\end{scope}
	\node[anchor=north east] at (1.5, 1.5) {$\Sigma_3$};
\end{tikzpicture}
&&\begin{tikzpicture}[scale=0.8]
	\clip (-0.03, -0.03) -- (-0.03, 2.97) -- (2.97, 2.97) -- (2.97, -0.03) -- cycle;
	\draw (0,0) -- (3, 0);
	\draw (0,0) -- (0, 3);
        \draw[very thick] (0,0) -- (6, 3);
        \draw (0,0) -- (3, 6);
	\begin{scope}\clip (0,0) -- (3, 0) -- (0, 3) -- cycle; \draw (0,0) circle (18pt);\end{scope}
        \begin{scope}\clip (0,0) -- (3, 0) -- (3, 6) -- cycle; \draw (0,0) circle (26pt);\end{scope}
	\node[anchor=north east] at (2.97, 2.97) {$\Theta_3$};
\end{tikzpicture}
&&
\begin{tikzpicture}[scale=0.8]
	\clip (-1.5, -1.5) -- (1.5, -1.5) -- (1.5, 1.5) -- (-1.5, 1.5) -- cycle;
	\fill[color=gray!30] (-3, 3) -- (3, -3) -- (-3, -3) -- cycle;
	\draw (1, 0) circle (2pt);
	\draw (0, 1) circle (2pt);
	\draw (0,0) -- (3,0);
	\draw (0,0) -- (0,3);
	\draw[very thick]  (0,0) -- (-4, 2);
	\draw[very thick]  (0,0) -- (2, -4);
        \fill (0,9pt) circle (2pt);
	\fill (9pt,0) circle (2pt);
        \fill (0,13pt) circle (2pt);
	\begin{scope}\clip (0,0) -- (1, -2) -- (2, 0) -- (0, 2) -- cycle; \draw (0,0) circle (9pt);\end{scope}
	\begin{scope}\clip (0,0) -- (0, 2) -- (-2, 1) -- cycle; \draw (0,0) circle (13pt);\end{scope}
	\begin{scope}\clip (0,0) -- (1, -2) -- (-2, -2) -- (-2, 1) -- cycle; \draw (0,0) circle (17pt);\end{scope}
	\node[anchor=north east] at (1.5, 1.5) {$\Sigma_4$};
\end{tikzpicture}
&&\begin{tikzpicture}[scale=0.8]
	\clip (-0.03, -0.03) -- (-0.03, 2.97) -- (2.97, 2.97) -- (2.97, -0.03) -- cycle;
	\draw (0,0) -- (3, 0);
	\draw (0,0) -- (0, 3);
        \draw (0,0) -- (6, 3);
        \draw[very thick] (0,0) -- (3, 6);
	\begin{scope}\clip (0,0) -- (3, 0) -- (0, 3) -- cycle; \draw (0,0) circle (18pt);\end{scope}
        \begin{scope}\clip (0,0) -- (6, 3) -- (0, 6) -- cycle; \draw (0,0) circle (26pt);\end{scope}
	\node[anchor=north east] at (2.97, 2.97) {$\Theta_4$};
\end{tikzpicture}
\\
\begin{tikzpicture}[scale=0.8]
	\clip (-1.5, -1.5) -- (1.5, -1.5) -- (1.5, 1.5) -- (-1.5, 1.5) -- cycle;
	\fill[color=gray!30] (-3, 3) -- (3, -3) -- (-3, -3) -- cycle;
	\draw (1, 0) circle (2pt);
	\draw (0, 1) circle (2pt);
	\draw (0,0) -- (3,0);
	\draw (0,0) -- (0,3);
	\draw (0,0)[very thick]  -- (-4, 2);
	\draw (0,0)[very thick]  -- (2, -4);
	\fill (9pt,0) circle (2pt);
        \fill (0,13pt) circle (2pt);
	\begin{scope}\clip (0,0) -- (1, -2) -- (2, 0) -- cycle; \draw (0,0) circle (9pt);\end{scope}
	\begin{scope}\clip (0,0) -- (0, 2) -- (-2, 1) -- cycle; \draw (0,0) circle (13pt);\end{scope}
	\begin{scope}\clip (0,0) -- (1, -2) -- (-2, -2) -- (-2, 1) -- cycle; \draw (0,0) circle (17pt);\end{scope}
	\node[anchor=north east] at (1.5, 1.5) {$\Sigma_5$};
\end{tikzpicture}
&&\begin{tikzpicture}[scale=0.8]
	\clip (-0.03, -0.03) -- (-0.03, 2.97) -- (2.97, 2.97) -- (2.97, -0.03) -- cycle;
	\draw (0,0) -- (3, 0);
	\draw (0,0) -- (0, 3);
        \draw (0,0) -- (6, 3);
        \draw (0,0) -- (3, 6);
	\begin{scope}\clip (0,0) -- (3, 0) -- (0, 3) -- cycle; \draw (0,0) circle (18pt);\end{scope}
	\begin{scope}\clip (0,0) -- (3, 0) -- (3, 6) -- cycle; \draw (0,0) circle (26pt);\end{scope}
        \begin{scope}\clip (0,0) -- (6, 3) -- (0, 6) -- cycle; \draw (0,0) circle (34pt);\end{scope}
	\node[anchor=north east] at (2.97, 2.97) {$\Theta_5$};
\end{tikzpicture}
\end{align*}
Moreover, we have the following colored fans which are not $\Pm$-fans.
According to Theorem~\ref{th:m}, these correspond to embeddings
which do not have the $A_2$-property.
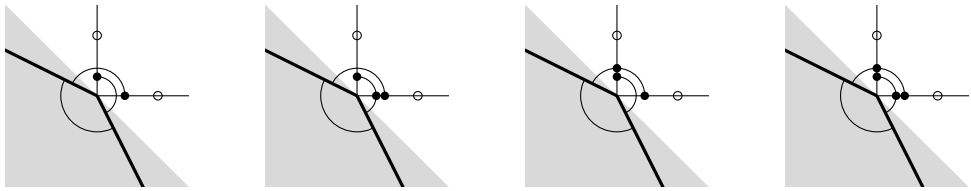
\begin{align*}
\begin{tikzpicture}[scale=0.8]
	\clip (-1.5, -1.5) -- (1.5, -1.5) -- (1.5, 1.5) -- (-1.5, 1.5) -- cycle;
	\fill[color=gray!30] (-3, 3) -- (3, -3) -- (-3, -3) -- cycle;
	\draw (1, 0) circle (2pt);
	\draw (0, 1) circle (2pt);
	\draw (0,0) -- (3,0);
	\draw (0,0) -- (0,3);
	\draw[very thick]  (0,0) -- (-4, 2);
	\draw[very thick]  (0,0) -- (2, -4);
	\fill (0,9pt) circle (2pt);
        \fill (13pt,0) circle (2pt);
	\begin{scope}\clip (0,0) -- (1, -2) -- (2, 0) -- (0, 2) -- cycle; \draw (0,0) circle (9pt);\end{scope}
	\begin{scope}\clip (0,0) -- (2, 0) -- (0, 2) -- (-2, 1) -- cycle; \draw (0,0) circle (13pt);\end{scope}
	\begin{scope}\clip (0,0) -- (1, -2) -- (-2, -2) -- (-2, 1) -- cycle; \draw (0,0) circle (17pt);\end{scope}
\end{tikzpicture}
&&
\begin{tikzpicture}[scale=0.8]
	\clip (-1.5, -1.5) -- (1.5, -1.5) -- (1.5, 1.5) -- (-1.5, 1.5) -- cycle;
	\fill[color=gray!30] (-3, 3) -- (3, -3) -- (-3, -3) -- cycle;
	\draw (1, 0) circle (2pt);
	\draw (0, 1) circle (2pt);
	\draw (0,0) -- (3,0);
	\draw (0,0) -- (0,3);
	\draw[very thick]  (0,0) -- (-4, 2);
	\draw[very thick]  (0,0) -- (2, -4);
	\fill (0,9pt) circle (2pt);
        \fill (13pt,0) circle (2pt);
	\fill (9pt,0) circle (2pt);
	\begin{scope}\clip (0,0) -- (1, -2) -- (2, 0) -- (0, 2) -- cycle; \draw (0,0) circle (9pt);\end{scope}
	\begin{scope}\clip (0,0) -- (2, 0) -- (0, 2) -- (-2, 1) -- cycle; \draw (0,0) circle (13pt);\end{scope}
	\begin{scope}\clip (0,0) -- (1, -2) -- (-2, -2) -- (-2, 1) -- cycle; \draw (0,0) circle (17pt);\end{scope}
\end{tikzpicture}
&&
\begin{tikzpicture}[scale=0.8]
	\clip (-1.5, -1.5) -- (1.5, -1.5) -- (1.5, 1.5) -- (-1.5, 1.5) -- cycle;
	\fill[color=gray!30] (-3, 3) -- (3, -3) -- (-3, -3) -- cycle;
	\draw (1, 0) circle (2pt);
	\draw (0, 1) circle (2pt);
	\draw (0,0) -- (3,0);
	\draw (0,0) -- (0,3);
	\draw[very thick]  (0,0) -- (-4, 2);
	\draw[very thick]  (0,0) -- (2, -4);
	\fill (0,9pt) circle (2pt);
        \fill (13pt,0) circle (2pt);
	\fill (0,13pt) circle (2pt);
	\begin{scope}\clip (0,0) -- (1, -2) -- (2, 0) -- (0, 2) -- cycle; \draw (0,0) circle (9pt);\end{scope}
	\begin{scope}\clip (0,0) -- (2, 0) -- (0, 2) -- (-2, 1) -- cycle; \draw (0,0) circle (13pt);\end{scope}
	\begin{scope}\clip (0,0) -- (1, -2) -- (-2, -2) -- (-2, 1) -- cycle; \draw (0,0) circle (17pt);\end{scope}
\end{tikzpicture}
&&
\begin{tikzpicture}[scale=0.8]
	\clip (-1.5, -1.5) -- (1.5, -1.5) -- (1.5, 1.5) -- (-1.5, 1.5) -- cycle;
	\fill[color=gray!30] (-3, 3) -- (3, -3) -- (-3, -3) -- cycle;
	\draw (1, 0) circle (2pt);
	\draw (0, 1) circle (2pt);
	\draw (0,0) -- (3,0);
	\draw (0,0) -- (0,3);
	\draw[very thick]  (0,0) -- (-4, 2);
	\draw[very thick]  (0,0) -- (2, -4);
	\fill (0,9pt) circle (2pt);
	\fill (9pt,0) circle (2pt);
        \fill (13pt,0) circle (2pt);
	\fill (0,13pt) circle (2pt);
	\begin{scope}\clip (0,0) -- (1, -2) -- (2, 0) -- (0, 2) -- cycle; \draw (0,0) circle (9pt);\end{scope}
	\begin{scope}\clip (0,0) -- (2, 0) -- (0, 2) -- (-2, 1) -- cycle; \draw (0,0) circle (13pt);\end{scope}
	\begin{scope}\clip (0,0) -- (1, -2) -- (-2, -2) -- (-2, 1) -- cycle; \draw (0,0) circle (17pt);\end{scope}
\end{tikzpicture}
\end{align*}
\end{example}

\begin{example}
\label{ex:2}
We consider the spherical homogeneous space
\begin{align*}
(\SL_2(\C)\times \C^*)/(T \times \{1\})
\end{align*}
where $T$ is a maximal torus in $\SL_2(\C)$. For details,
we refer, for instance, to \cite[Example~2.5.1]{pez10}.
The following pictures show the vector spaces $\Nm_\Q$
and $K_\Q$ (with the same notation as in Example~\ref{ex:1}). For
the two colors, $\Dm = \{D_1, D_2\}$, we have $\rho(D_1) = \rho(D_2)$,
so that in order to be able to distinguish them, we have put one circle together with the number $1$ below and one circle together with the number $2$ above the ray.
We have added two $\SL_2(\C) \times \C^*$-invariant prime divisors,
$\gdiv = \{D_3, D_4\}$. Note that we have $\Qm(D_3) = \Qm(D_4)$, see the paragraph before Remark~\ref{rem:sphcox}.
\begin{align*}
\begin{tikzpicture}[scale=0.8]
	\clip (-1.5, -1.5) -- (1.5, -1.5) -- (1.5, 1.5) -- (-1.5, 1.5) -- cycle;
	\fill[color=gray!30] (0, -3) -- (0, 3) -- (-3, 3) -- (-3, -3) -- cycle;
	\draw (1, -0.07) circle (2pt);
	\draw (1, 0.07) circle (2pt);
	\draw (0,0) -- (1.5,0) node[very near end,below]{${\scriptstyle 1}$} node[very near end,above]{${\scriptstyle 2}$};
	\draw (0,0) -- (-1.5,1.5) node[very near end,right]{${\scriptstyle 3}$};
	\draw (0,0) -- (-1.5, -1.5) node[very near end,right]{${\scriptstyle 4}$};;
\end{tikzpicture}
&&
\begin{tikzpicture}[scale=0.8]
	\clip (-0.03, -0.03) -- (-0.03, 2.97) -- (2.97, 2.97) -- (2.97, -0.03) -- cycle;
	\draw (0,0) -- (3, 0) node[very near end,above]{${\scriptstyle 2}$};
	\draw (0,0) -- (0, 3) node[very near end,right]{${\scriptstyle 1}$};
        \draw (0,0) -- (3, 3) node[very near end,left]{${\scriptstyle 3,4}$};
\end{tikzpicture}
\end{align*}
The following pictures show all possible true maximal $\Pm$-fans $\Sigma_j$ and their corresponding true maximal $\Qm$-bunches $\Theta_j$.
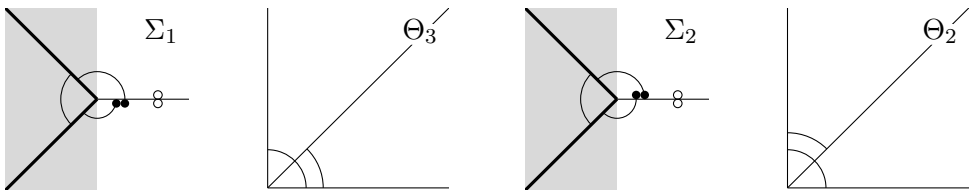
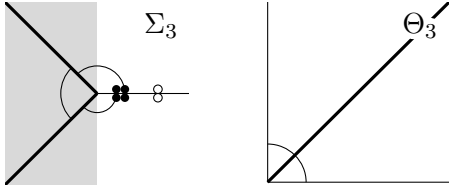
\begin{align*}
\begin{tikzpicture}[scale=0.8]
	\clip (-1.5, -1.5) -- (1.5, -1.5) -- (1.5, 1.5) -- (-1.5, 1.5) -- cycle;
	\fill[color=gray!30] (0, -3) -- (0, 3) -- (-3, 3) -- (-3, -3) -- cycle;
	\draw (1, -0.07) circle (2pt);
	\draw (1, 0.07) circle (2pt);
	\draw (0,0) -- (3,0);
	\draw[very thick] (0,0) -- (-3,3);
	\draw[very thick] (0,0) -- (-3, -3);
	\fill (9pt,-0.07) circle (2pt);
        \fill (13pt,-0.07) circle (2pt);
	\begin{scope}\clip (0,0) -- (-3, -3) -- (3, 0) -- cycle; \draw (0,0) circle (9pt);\end{scope}
	\begin{scope}\clip (0,0) -- (3, 0) -- (-3, 3) -- cycle; \draw (0,0) circle (13pt);\end{scope}
	\begin{scope}\clip (0,0) -- (-3, 3) -- (-3, -3) -- cycle; \draw (0,0) circle (17pt);\end{scope}
	\node[anchor=north east] at (1.5, 1.5) {$\Sigma_1$};
\end{tikzpicture}
&&
\begin{tikzpicture}[scale=0.8]
	\clip (-0.03, -0.03) -- (-0.03, 2.97) -- (2.97, 2.97) -- (2.97, -0.03) -- cycle;
	\draw (0,0) -- (3, 0);
	\draw (0,0) -- (0, 3);
        \draw (0,0) -- (6, 6);
	\begin{scope}\clip (0,0) -- (3, 0) -- (0, 3) -- cycle; \draw (0,0) circle (18pt);\end{scope}
	\begin{scope}\clip (0,0) -- (3, 0) -- (6, 6) -- cycle; \draw (0,0) circle (26pt);\end{scope}
	\node[anchor=north east] at (2.97, 2.97) {\contour{white}{$\Theta_3$}};
\end{tikzpicture}
&&
\begin{tikzpicture}[scale=0.8]
	\clip (-1.5, -1.5) -- (1.5, -1.5) -- (1.5, 1.5) -- (-1.5, 1.5) -- cycle;
	\fill[color=gray!30] (0, -3) -- (0, 3) -- (-3, 3) -- (-3, -3) -- cycle;
	\draw (1, -0.07) circle (2pt);
	\draw (1, 0.07) circle (2pt);
	\draw (0,0) -- (3,0);
	\draw[very thick]  (0,0) -- (-3,3);
	\draw[very thick]  (0,0) -- (-3, -3);
	\fill (9pt,0.07) circle (2pt);
        \fill (13pt,0.07) circle (2pt);
	\begin{scope}\clip (0,0) -- (-3, -3) -- (3, 0) -- cycle; \draw (0,0) circle (9pt);\end{scope}
	\begin{scope}\clip (0,0) -- (3, 0) -- (-3, 3) -- cycle; \draw (0,0) circle (13pt);\end{scope}
	\begin{scope}\clip (0,0) -- (-3, 3) -- (-3, -3) -- cycle; \draw (0,0) circle (17pt);\end{scope}
	\node[anchor=north east] at (1.5, 1.5) {$\Sigma_2$};
\end{tikzpicture}
&&
\begin{tikzpicture}[scale=0.8]
	\clip (-0.03, -0.03) -- (-0.03, 2.97) -- (2.97, 2.97) -- (2.97, -0.03) -- cycle;
	\draw (0,0) -- (3, 0);
	\draw (0,0) -- (0, 3);
        \draw (0,0) -- (6, 6);
	\begin{scope}\clip (0,0) -- (3, 0) -- (0, 3) -- cycle; \draw (0,0) circle (18pt);\end{scope}
	\begin{scope}\clip (0,0) -- (6, 6) -- (0, 3) -- cycle; \draw (0,0) circle (26pt);\end{scope}
	\node[anchor=north east] at (2.97, 2.97) {\contour{white}{$\Theta_2$}};
\end{tikzpicture}
\\
\begin{tikzpicture}[scale=0.8]
	\clip (-1.5, -1.5) -- (1.5, -1.5) -- (1.5, 1.5) -- (-1.5, 1.5) -- cycle;
	\fill[color=gray!30] (0, -3) -- (0, 3) -- (-3, 3) -- (-3, -3) -- cycle;
	\draw (1, -0.07) circle (2pt);
	\draw (1, 0.07) circle (2pt);
	\draw (0,0) -- (3,0);
	\draw[very thick]  (0,0) -- (-3,3);
	\draw[very thick]  (0,0) -- (-3, -3);
	\fill (9pt,-0.07) circle (2pt);
        \fill (13pt,-0.07) circle (2pt);
	\fill (9pt,0.07) circle (2pt);
        \fill (13pt,0.07) circle (2pt);
	\begin{scope}\clip (0,0) -- (-3, -3) -- (3, 0) -- cycle; \draw (0,0) circle (9pt);\end{scope}
	\begin{scope}\clip (0,0) -- (3, 0) -- (-3, 3) -- cycle; \draw (0,0) circle (13pt);\end{scope}
	\begin{scope}\clip (0,0) -- (-3, 3) -- (-3, -3) -- cycle; \draw (0,0) circle (17pt);\end{scope}
	\node[anchor=north east] at (1.5, 1.5) {$\Sigma_3$};
\end{tikzpicture}
&&
\begin{tikzpicture}[scale=0.8]
	\clip (-0.03, -0.03) -- (-0.03, 2.97) -- (2.97, 2.97) -- (2.97, -0.03) -- cycle;
	\draw (0,0) -- (3, 0);
	\draw (0,0) -- (0, 3);
        \draw[very thick] (0,0) -- (6, 6);
	\begin{scope}\clip (0,0) -- (3, 0) -- (0, 3) -- cycle; \draw (0,0) circle (18pt);\end{scope}
	\node[anchor=north east] at (2.97, 2.97) {\contour{white}{$\Theta_3$}};
\end{tikzpicture}
\end{align*}
Moreover, we have the following colored fans which are not $\Pm$-fans.
According to Theorem~\ref{th:m}, these correspond to embeddings
which do not have the $A_2$-property.
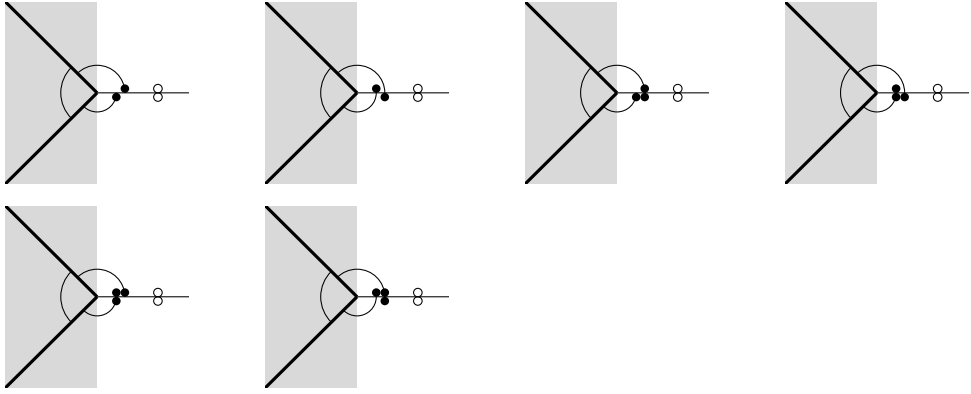
\begin{align*}
\begin{tikzpicture}[scale=0.8]
	\clip (-1.5, -1.5) -- (1.5, -1.5) -- (1.5, 1.5) -- (-1.5, 1.5) -- cycle;
	\fill[color=gray!30] (0, -3) -- (0, 3) -- (-3, 3) -- (-3, -3) -- cycle;
	\draw (1, -0.07) circle (2pt);
	\draw (1, 0.07) circle (2pt);
	\draw (0,0) -- (3,0);
	\draw[very thick] (0,0) -- (-3,3);
	\draw[very thick]  (0,0) -- (-3, -3);
	\fill (9pt,-0.07) circle (2pt);
        \fill (13pt,0.07) circle (2pt);
	\begin{scope}\clip (0,0) -- (-3, -3) -- (3, 0) -- cycle; \draw (0,0) circle (9pt);\end{scope}
	\begin{scope}\clip (0,0) -- (3, 0) -- (-3, 3) -- cycle; \draw (0,0) circle (13pt);\end{scope}
	\begin{scope}\clip (0,0) -- (-3, 3) -- (-3, -3) -- cycle; \draw (0,0) circle (17pt);\end{scope}
\end{tikzpicture}
&&
\begin{tikzpicture}[scale=0.8]
	\clip (-1.5, -1.5) -- (1.5, -1.5) -- (1.5, 1.5) -- (-1.5, 1.5) -- cycle;
	\fill[color=gray!30] (0, -3) -- (0, 3) -- (-3, 3) -- (-3, -3) -- cycle;
	\draw (1, -0.07) circle (2pt);
	\draw (1, 0.07) circle (2pt);
	\draw (0,0) -- (3,0);
	\draw[very thick]  (0,0) -- (-3,3);
	\draw[very thick]  (0,0) -- (-3, -3);
        \fill (13pt,-0.07) circle (2pt);
	\fill (9pt,0.07) circle (2pt);
	\begin{scope}\clip (0,0) -- (-3, -3) -- (3, 0) -- cycle; \draw (0,0) circle (9pt);\end{scope}
	\begin{scope}\clip (0,0) -- (3, 0) -- (-3, 3) -- cycle; \draw (0,0) circle (13pt);\end{scope}
	\begin{scope}\clip (0,0) -- (-3, 3) -- (-3, -3) -- cycle; \draw (0,0) circle (17pt);\end{scope}
\end{tikzpicture}
&&
\begin{tikzpicture}[scale=0.8]
	\clip (-1.5, -1.5) -- (1.5, -1.5) -- (1.5, 1.5) -- (-1.5, 1.5) -- cycle;
	\fill[color=gray!30] (0, -3) -- (0, 3) -- (-3, 3) -- (-3, -3) -- cycle;
	\draw (1, -0.07) circle (2pt);
	\draw (1, 0.07) circle (2pt);
	\draw (0,0) -- (3,0);
	\draw[very thick]  (0,0) -- (-3,3);
	\draw[very thick]  (0,0) -- (-3, -3);
	\fill (9pt,-0.07) circle (2pt);
        \fill (13pt,-0.07) circle (2pt);
        \fill (13pt,0.07) circle (2pt);
	\begin{scope}\clip (0,0) -- (-3, -3) -- (3, 0) -- cycle; \draw (0,0) circle (9pt);\end{scope}
	\begin{scope}\clip (0,0) -- (3, 0) -- (-3, 3) -- cycle; \draw (0,0) circle (13pt);\end{scope}
	\begin{scope}\clip (0,0) -- (-3, 3) -- (-3, -3) -- cycle; \draw (0,0) circle (17pt);\end{scope}
\end{tikzpicture}
&&
\begin{tikzpicture}[scale=0.8]
	\clip (-1.5, -1.5) -- (1.5, -1.5) -- (1.5, 1.5) -- (-1.5, 1.5) -- cycle;
	\fill[color=gray!30] (0, -3) -- (0, 3) -- (-3, 3) -- (-3, -3) -- cycle;
	\draw (1, -0.07) circle (2pt);
	\draw (1, 0.07) circle (2pt);
	\draw (0,0) -- (3,0);
	\draw[very thick]  (0,0) -- (-3,3);
	\draw[very thick]  (0,0) -- (-3, -3);
	\fill (9pt,-0.07) circle (2pt);
        \fill (13pt,-0.07) circle (2pt);
	\fill (9pt,0.07) circle (2pt);
	\begin{scope}\clip (0,0) -- (-3, -3) -- (3, 0) -- cycle; \draw (0,0) circle (9pt);\end{scope}
	\begin{scope}\clip (0,0) -- (3, 0) -- (-3, 3) -- cycle; \draw (0,0) circle (13pt);\end{scope}
	\begin{scope}\clip (0,0) -- (-3, 3) -- (-3, -3) -- cycle; \draw (0,0) circle (17pt);\end{scope}
\end{tikzpicture}
\\
\begin{tikzpicture}[scale=0.8]
	\clip (-1.5, -1.5) -- (1.5, -1.5) -- (1.5, 1.5) -- (-1.5, 1.5) -- cycle;
	\fill[color=gray!30] (0, -3) -- (0, 3) -- (-3, 3) -- (-3, -3) -- cycle;
	\draw (1, -0.07) circle (2pt);
	\draw (1, 0.07) circle (2pt);
	\draw (0,0) -- (3,0);
	\draw[very thick]  (0,0) -- (-3,3);
	\draw[very thick]  (0,0) -- (-3, -3);
	\fill (9pt,-0.07) circle (2pt);
	\fill (9pt,0.07) circle (2pt);
        \fill (13pt,0.07) circle (2pt);
	\begin{scope}\clip (0,0) -- (-3, -3) -- (3, 0) -- cycle; \draw (0,0) circle (9pt);\end{scope}
	\begin{scope}\clip (0,0) -- (3, 0) -- (-3, 3) -- cycle; \draw (0,0) circle (13pt);\end{scope}
	\begin{scope}\clip (0,0) -- (-3, 3) -- (-3, -3) -- cycle; \draw (0,0) circle (17pt);\end{scope}
\end{tikzpicture}
&&
\begin{tikzpicture}[scale=0.8]
	\clip (-1.5, -1.5) -- (1.5, -1.5) -- (1.5, 1.5) -- (-1.5, 1.5) -- cycle;
	\fill[color=gray!30] (0, -3) -- (0, 3) -- (-3, 3) -- (-3, -3) -- cycle;
	\draw (1, -0.07) circle (2pt);
	\draw (1, 0.07) circle (2pt);
	\draw (0,0) -- (3,0);
	\draw[very thick]  (0,0) -- (-3,3);
	\draw[very thick]  (0,0) -- (-3, -3);
        \fill (13pt,-0.07) circle (2pt);
	\fill (9pt,0.07) circle (2pt);
        \fill (13pt,0.07) circle (2pt);
	\begin{scope}\clip (0,0) -- (-3, -3) -- (3, 0) -- cycle; \draw (0,0) circle (9pt);\end{scope}
	\begin{scope}\clip (0,0) -- (3, 0) -- (-3, 3) -- cycle; \draw (0,0) circle (13pt);\end{scope}
	\begin{scope}\clip (0,0) -- (-3, 3) -- (-3, -3) -- cycle; \draw (0,0) circle (17pt);\end{scope}
\end{tikzpicture}
\end{align*}
\end{example}

\section{Proof of Theorem~\ref{thm:fb}}
\label{sec:p}
This section is a shortened and suitably modified version of \cite[2.2.3]{coxrings}.
We define
\begin{align*}
\delta \coloneqq \cone(e^*_D : D \in \adiv)\text{,} &&
\gamma \coloneqq \cone(e_D : D \in \adiv)\text{.}
\end{align*}
These cones are dual to each other, and we have the face correspondence,
\ie mutually inverse bijections
\begin{align*}
\operatorname{faces}(\delta) &\leftrightarrow \operatorname{faces}(\gamma)\text{,}\\
\delta_0 &\mapsto \delta_0^* \coloneqq \delta_0^\perp \cap \gamma\text{,}\\
\gamma_0^\perp \cap \delta \eqqcolon \gamma_0^* &\mapsfrom \gamma_0\text{.}
\end{align*}

\begin{definition}
A face $\delta_0 \preceq \delta$
is called \emph{supported} if
$\Pm(\delta_0)^\circ \cap \Vm \ne \emptyset$.
A face $\gamma_0 \preceq \gamma$
is called \emph{supported} if $\gamma_0^* \preceq \delta$ is supported.
A \emph{$\delta$-collection} (resp.~a \emph{$\gamma$-collection})
is a set of supported faces of $\delta$ (resp.~of $\gamma$).
A \emph{$\Pm$-collection} (resp.~a \emph{$\Qm$-collection}) is
a set of supported $\Pm$-cones (resp.~of supported $\Qm$-cones).
\end{definition}

\begin{remark}
A $\Qm$-cone $\tau$ is supported if and only if there exists a
supported face $\gamma_0 \preceq \gamma$ with $\Qm(\gamma_0) = \tau$.
\end{remark}

The idea of the proof is to decompose the $\sharp$-operation between
$\Qm$-collections and $\Pm$-collections according to the following
scheme of further operations.
\begin{align*}
\xymatrix{
\{\text{$\gamma$-collections}\} \ar@{<->}[r]^-{*} \ar@{<->}[d]^-{\Qm_{\smash{\downarrow}}}_-{\Qm^{\smash{\uparrow}}} & \{\text{$\delta$-collections}\} 
\ar@{<->}[d]_-{\Pm_{\smash{\downarrow}}}^-{\Pm^{\smash{\uparrow}}} 
\\
\{\text{$\Qm$-collections}\} \ar@{<->}[r]_-{\sharp} & \{\text{$\Pm$-collections}\}
}
\end{align*}

\begin{definition}
An \emph{$L_\Q$-invariant separating linear form} for two faces $\delta_1, \delta_2 \preceq \delta$
is an element $e \in \Q^{\adiv}$ such that
\begin{align*}
e|_{L_\Q} = 0\text{,} && e|_{\delta_1} \ge 0\text{,} && e|_{\delta_2} \le 0\text{,} &&
\delta_1 \cap e^\perp = \delta_1 \cap \delta_2 = e^\perp \cap \delta_2\text{.}
\end{align*}
\end{definition}

\begin{definition}
A $\delta$-collection $\Af$ is called
\begin{enumerate}[(i)]
\item \emph{separated} if any two $\delta_1, \delta_2 \in \Af$ admit an $L_\Q$-invariant separating linear form,
\item \emph{saturated} if for any $\delta_1 \in \Af$ every supported $\delta_2 \preceq \delta_1$ which 
is $L_\Q$-invariantly separable from $\delta_1$ also belongs to $\Af$,
\item \emph{true} if we have $0 \in \Af$ for $\Dm \ne \emptyset$ and $\cone(e^*_D) \in \Af$ for every $D \in \gdiv$,
\item \emph{maximal} if it is maximal among the separated $\delta$-collections.
\end{enumerate}
\end{definition}

\begin{definition}
A $\gamma$-collection $\Bf$ is called
\begin{enumerate}[(i)]
\item \emph{connected} if for any $\gamma_1, \gamma_2 \in \Bf$ we have
$\Qm(\gamma_1)^\circ \cap \Qm(\gamma_2)^\circ \ne \emptyset$,
\item \emph{saturated} if for any $\gamma_1 \in \Bf$ every supported $\gamma_2$
with $\gamma \succeq \gamma_2 \succeq \gamma_1$ and $\Qm(\gamma_1)^\circ \subseteq \Qm(\gamma_2)^\circ$
also belongs to $\Bf$,
\item \emph{true} if we have $\gamma \in \Bf$ for $\Dm \ne \emptyset$ and $\cone(e^*_D)^* \in \Bf$
for every $D \in \gdiv$,
\item \emph{maximal} if it is maximal among the connected $\gamma$-collections.
\end{enumerate}
\end{definition}

\begin{prop}
\label{prop:mid}
We have mutually inverse bijections sending separated (saturated, true, maximal)
collections to connected (saturated, true, maximal) collections:
\begin{align*}
\{\text{separated $\delta$-collections}\} &\leftrightarrow \{\text{connected $\gamma$-collections}\}\text{,}\\
\Af & \mapsto \Af^* \coloneqq \{\delta_0^* : \delta_0 \in \Af\}\text{,}\\
\{\gamma_0^* : \gamma_0 \in \Bf\} \eqqcolon \Bf^* &\mapsfrom \Bf\text{.}
\end{align*}
\end{prop}
\begin{proof}
See \cite[Proposition~2.2.3.5]{coxrings}.
\end{proof}

\begin{definition}
A $\Qm$-collection $\Theta$ is called
\begin{enumerate}[(i)]
\item \emph{connected} if for any $\tau_1, \tau_2 \in \Theta$ we have
$\tau_1^\circ \cap \tau_2^\circ \ne \emptyset$,
\item \emph{saturated} if for any $\tau \in \Theta$ every supported $\Qm$-cone $\tau_0$
with $\tau^\circ \subseteq \tau_0^\circ$ also belongs to $\Theta$,
\item \emph{true} if we have $\cone(\Qm(\adiv)) \in \Theta$ for $\Dm \ne \emptyset$ and 
$\cone(\Qm(\adiv \setminus \{D\})) \in \Theta$ for every $D \in \gdiv$,
\item \emph{maximal} if it is maximal among the connected $\Qm$-collections.
\end{enumerate}
\end{definition}

\begin{definition}
We define the $\Qm$-lift and the $\Qm$-drop to be the maps
\begin{align*}
\Qm^\uparrow\colon \{\text{$\Qm$-collections}\} &\to \{\text{$\gamma$-collections}\}\text{,}\\
\Theta &\mapsto \Qm^\uparrow\Theta \coloneqq \{\gamma_0 \preceq \gamma : \text{$\gamma_0$ is supported and
$\Qm(\gamma_0) \in \Theta$}\}\text{,}\\
\Qm_{\downarrow}\colon \{\text{$\gamma$-collections}\} &\to \{\text{$\Qm$-collections}\}\text{,}\\
\Bf &\mapsto \Qm_{\downarrow}\Bf \coloneqq \{\Qm(\gamma_0) : \gamma_0 \in \Bf\}\text{.}
\end{align*}
\end{definition}

\begin{prop}
\label{prop:Qlift}
The $\Qm$-lift is injective and sends connected (saturated, true, maximal) $\Qm$-collections
to connected (saturated, true, maximal) $\gamma$-collections. Moreover, we have mutually
inverse bijections sending true collections to true collections:
\begin{align*}
\{\text{maximal $\Qm$-collections}\} &\leftrightarrow \{\text{maximal $\gamma$-collections}\}\text{,}\\
\Theta &\mapsto \Qm^{\uparrow}\Theta\text{,}\\
\Qm_{\downarrow}\Bf &\mapsfrom \Bf\text{.}
\end{align*}
\end{prop}
\begin{proof}
See \cite[Proposition~2.2.3.8]{coxrings}.
\end{proof}

\begin{definition}
A $\Pm$-collection $\Sigma$ is called
\begin{enumerate}[(i)]
\item \emph{separated} if any two $(\Cm_1, \Fm_1), (\Cm_2, \Fm_2)\in \Sigma$
intersect in a common face,
\item \emph{saturated} if for any $(\Cm, \Fm) \in \Sigma$ every supported
face of $(\Cm, \Fm)$ also is in $\Sigma$,
\item \emph{true} if we have $(0, \emptyset) \in \Sigma$ for $\bdiv \ne \emptyset$ and
$(\cone(\Pm(D)), \emptyset) \in \Sigma$ for every $D \in \gdiv$,
\item \emph{maximal} if it is maximal among the separated $\Pm$-collections.
\end{enumerate}
\end{definition}

\begin{definition}
We define the $\Pm$-lift and the $\Pm$-drop to be the maps
\begin{align*}
\Pm^\uparrow\colon \{\text{$\Pm$-collections}\} &\to \{\text{$\delta$-collections}\}\text{,}\\
\Sigma &\mapsto \Pm^\uparrow\Sigma \coloneqq \{\delta_0 \preceq \delta : \text{$(\Pm(\delta_0),\{D \in \bdiv : e^*_D \in \delta_0\}) \in \Sigma$}\}\text{,}\\
\Pm_{\downarrow}\colon \{\text{$\delta$-collections}\} &\to \{\text{$\Pm$-collections}\}\text{,}\\
\Af &\mapsto \Pm_{\downarrow}\Af \coloneqq \{(\Pm(\delta_0), \{D \in \bdiv : e^*_D \in \delta_0\}) : \delta_0 \in \Af\}\text{.}
\end{align*}
\end{definition}

\begin{prop}
\label{prop:Pdrop}
The $\Pm$-drop is surjective and sends separated (saturated, true, maximal)
$\delta$-collections to separated (saturated, true, maximal) $\Pm$-collections.
If the elements in $\Pm(\gdiv)$ generate pairwise different rays,
then we have mutually inverse
bijections sending saturated (maximal) collections
to saturated (maximal) collections:
\begin{align*}
\{\text{true separated $\delta$-collections}\} &\leftrightarrow \{\text{true separated $\Pm$-collections}\}\text{,}\\
\Af &\mapsto \Pm_{\downarrow}\Af\text{,}\\
\Pm^{\uparrow}\Sigma &\mapsfrom \Sigma\text{.}
\end{align*}
\end{prop}
\begin{proof}
For every $\Pm$-collection $\Sigma$ we have $\Sigma = \Pm_{\downarrow}\Pm^{\uparrow}\Sigma$.
In particular, $\Pm_{\downarrow}$ is surjective.
Let $\delta_1, \delta_2 \preceq \delta$ be two faces
admitting an $L_\Q$-invariant separating linear form $e \in \Q^{\adiv}$,
and define
\begin{align*}
(\Cm_1, \Fm_1) &\coloneqq (\Pm(\delta_1),\{D \in \Dm : e^*_D \in \delta_1\})\text{,}\\
(\Cm_2, \Fm_2) &\coloneqq (\Pm(\delta_2),\{D \in \Dm : e^*_D \in \delta_2\})\text{.}
\end{align*}
Then $e$ can be interpreted as $e \in \Mm_\Q$ with
$e|_{\Cm_1} \ge 0$, $e|_{\Cm_2} \le 0$,
\begin{align*}
\Cm_1 \cap e^\perp = \Cm_1 \cap \Cm_2 = e^\perp \cap \Cm_2\text{,} &&\text{and}&&
\Fm_1 \cap \rho|_{\Dm}^{-1}(e^\perp) = \rho|_{\Dm}^{-1}(e^\perp) \cap \Fm_2\text{.}
\end{align*}
It now follows from Remark~\ref{rem:intf} that $\Pm_{\downarrow}$ preserves separatedness and saturatedness.
The fact that $\Pm_{\downarrow}$ preserves the properties true and maximal is obvious.

Now assume that the elements in $\Pm(\gdiv)$ generate
pairwise different rays. Consider
a true separated $\Pm$-collection $\Sigma$. Then, for every
$(\Cm, \Fm) \in \Sigma$ and every $D \in \gdiv$ we have $\Pm(D) \in \Cm$ if and only if
$\Q_{\ge 0}\cdot \Pm(D)$ is an extremal ray of $\Cm$.
Moreover, for every $D' \in \bdiv$ such that there 
exists $D'' \in \gdiv$ with $\Q_{\ge 0}\cdot \Pm(D') = \Q_{\ge 0}\cdot \Pm(D'')$
we have $D' \notin \Fm$.
Consequently, for every
$(\Cm, \Fm) \in \Sigma$ there is a unique $\delta_0 \preceq \delta$ with
$(\Pm(\delta_0), \{D \in \Dm : e^*_D \in \delta_0\}) = (\Cm, \Fm)$.
It follows that $\Pm^{\uparrow}\Sigma$ is true and separated, and, if $\Sigma$
is saturated (maximal), then $\Pm^{\uparrow}\Sigma$ is also saturated (maximal).
Moreover, we conclude that $\Pm_{\downarrow}$ restricted to the true separated
$\delta$-collections is injective.
\end{proof}

\begin{proof}[Proof of Theorem~\ref{thm:fb}]
First, observe that the (true, maximal) $\Qm$-bunches
are precisely the (true, maximal) connected saturated
$\Qm$-collections and the (true, maximal) $\Pm$-quasifans
are precisely the (true, maximal) separated saturated $\Pm$-collections.
Next observe that we have
\begin{align*}
\Theta^\sharp = \Pm_{\downarrow}((\Qm^{\uparrow}\Theta)^*)\text{,}
&& \Sigma^\sharp = \Qm_{\downarrow}((\Pm^{\uparrow}\Sigma)^*)\text{.}
\end{align*}
The first part of Theorem~\ref{thm:fb} now follows
from Propositions~\ref{prop:Qlift}, \ref{prop:mid}, and \ref{prop:Pdrop}.

It remains to show that $\Pm$-fans consisting of simplicial $\Pm$-cones
correspond to $\Qm$-bunches consisting of full-dimensional $\Qm$-cones.
A true $\Pm$-fan $\Sigma$ is simplicial exactly when for every
$(\cone(I), I \cap \bdiv) \in \Sigma$ with $I \subseteq \adiv$
and any subset $I_0 \subseteq I$ we have that
$(\cone(I_0), I_0 \cap \bdiv)$ is a face of $(\cone(I), I \cap \bdiv)$.
This means that for every $\tau = \cone(\Qm(J))$ with $J \subseteq \adiv$ in
the corresponding true $\Qm$-bunch $\Theta$ and every $J_1$ with $J \subseteq J_1 \subseteq \adiv$
we have $\tau^\circ \subseteq \cone(\Qm(J_1))^\circ$. Because
the vectors $\{\Qm(D) : D \in \adiv\}$ generate $K_\Q$, this is exactly
the case when all cones in $\Theta$ are of full dimension.
\end{proof}

\section{Bunched rings}
\label{sec:br}

Let $X_0$ be a normal irreducible variety with 
$\Gamma(X_0, \Om^*_{X_0})=\C^*$,
finitely generated divisor class group $K \coloneqq \Cl(X_0)$,
and finitely generated \emph{Cox ring} 
\begin{align*}
R \coloneqq \Rm(X_0) \coloneqq \bigoplus_{[D]\in K} \Gamma(X_0, \Om_{X_0}(D))\text{,}
\end{align*}
where some care has to be taken in order to define the multiplication
law. The Cox ring $R$ is factorially $K$-graded.
This means that every homogeneous nonzero nonunit in $R$ can be written
as a product of $K$-primes, where a
$K$-prime is a homogeneous nonzero nonunit
$f \in R$ such that $f \mid gh$ with homogeneous $g,h \in R$
always implies $f \mid g$ or $f \mid h$.
For details, we refer to \cite[1.4 and 1.6]{coxrings}.

With $\overline{X} \coloneqq \Spec R$ the $K$-grading
on $R$ corresponds to an $S$-action on $\overline{X}$
where $S \coloneqq \Spec \C[K]$ is a quasitorus (\ie a diagonalizable group) with character group $K$.
There exists an open $S$-stable subvariety $\smash{\widehat{X}}_0 \subseteq \overline{X}$
with complement of codimension at least $2$ such that we obtain
a good quotient
$\pi\colon \smash{\widehat{X}}_0 \to X_0$ for the $S$-action.

If $X$ is any other normal irreducible variety with the same graded Cox ring as $X_0$,
we obtain $X$ as a good quotient $\pi\colon \smash{\widehat{X}} \to X$
for some other open subvariety $\smash{\widehat{X}} \subseteq \overline{X}$ with complement
of codimension at least $2$ (the varieties $X_0$ and $X$ also differ
at most in codimension $2$). The theory of bunched rings
(which first appeared in \cite{cc1, cc2}) can be used
to find such $\smash{\widehat{X}} \subseteq \overline{X}$ provided that the quotient $X$
has the $A_2$-property.
We recall some definitions and results on bunched rings from
\cite[3.2]{coxrings}.

\begin{definition}[{\cite[Definitions~3.2.1.1 and 3.2.1.2]{coxrings}}]
\label{def:prb}
Let $\Ff$ be a finite system of pairwise
nonassociated $K$-prime generators for $R$.
\begin{enumerate}[(i)]
\item The $K$-grading is said to be \emph{almost free} if
for every $f_0 \in \Ff$ the set \begin{align*}\{\deg f : f \in \Ff \setminus \{f_0\}\}\end{align*}
generates $K$ as an abelian group.
\item The set of \emph{projected $\Ff$-faces} is the set
\begin{align*}
\Omega_{\Ff} \coloneqq \rleft\{\cone(\deg f : f \in J) : J \subseteq \Ff, \bigcap_{f \notin J} \V(f) \mathbin{\Big\backslash} \bigcup_{f\in J} \V(f) \ne \emptyset\rright\}
\end{align*}
of cones in $K_\Q$ where the $\V(f)$ are considered as subsets of $\overline{X}$.
\item An \emph{$\Ff$-bunch} is a nonempty subset $\Theta \subseteq \Omega_{\Ff}$ such that
\begin{enumerate}
\item for any $\tau_1, \tau_2 \in \Theta$ we have $\tau_1^\circ \cap \tau_2^\circ = \emptyset$,
\item for any $\tau \in \Theta$ every $\tau_0 \in \Omega_{\Ff}$ with $\tau^\circ \subseteq \tau_0^\circ$
also belongs to $\Theta$.
\end{enumerate}
\item An $\Ff$-bunch $\Theta$ is called \emph{true} if it contains
\begin{align*}\cone(\deg f : f \in \Ff \setminus \{f_0\})\end{align*}
for every $f_0 \in \Ff$.
\end{enumerate}
If $\Theta$ is an $\Ff$-bunch, the triple $(R, \Ff, \Theta)$ is called a \emph{bunched ring}.
\end{definition}

As in \cite[3.2.1]{coxrings}, for 
any $\Ff$-bunch $\Theta$ we set
\begin{align*}
\smash{\widehat{X}} \coloneqq \bigcup_{\tau \in \Theta}
\bigcup_{\substack{J \subseteq \Ff\\\cone(\deg f:f\in J)=\tau}}
\rleft(\overline{X}\setminus \bigcup_{f\in J}\V(f)\rright)\text{.}
\end{align*}
Then there exists a good quotient $\pi\colon \widehat{X} \to X$ and we say
that the variety $X$ arises from the bunched ring $(R, \Ff, \Theta)$.
As in \cite[3.3.1]{coxrings}, for any $\tau \in \Omega_{\Ff}$ we define
\begin{align*}
X(\tau) \coloneqq \bigcup_{\substack{J \subseteq \Ff\\\cone(\deg f:f\in J)=\tau}}\rleft(\bigcap_{f \notin J} \pi(\V(f)) \mathbin{\Big\backslash} \bigcup_{f\in J} \pi(\V(f))\rright)\subseteq X
\end{align*}
where the $\V(f)$ are now considered as subsets of $\smash{\widehat{X}}$.
We obtain a disjoint union
\begin{align*}
X = \bigcup_{\tau \in \Theta} X(\tau)\text{.}
\end{align*}
 Note that
we have $X(\tau) \ne \emptyset$ if and only if $\tau \in \Theta$.

\begin{prop}[{\cite[Proposition~3.2.1.9]{coxrings}}]
\label{prop:ea2}
Let $X$ be a variety with graded Cox ring $R$.
Then $X$ has the $A_2$-property
if and only if it is
an open subvariety of a variety
arising from a bunched ring $(R, \Ff, \Theta)$.
\end{prop}

This can be generalized to take into
account the $A_k$-property for any $k \ge 2$ in a straightforward manner.

\begin{prop}[{\cite[Exercise~3.5(4)]{coxrings}}]
\label{prop:eak3}
Let $X$ be a variety with graded Cox ring $R$,
and let $k \ge 2$.
Then $X$ has the $A_k$-property
if and only if it is
an open subvariety of a variety
arising from a bunched ring $(R, \Ff, \Theta)$
such that for any $k$ cones
$\tau_1, \dots, \tau_k \in \Theta$ we have $\tau_1^\circ \cap \dots \cap \tau_k^\circ \ne \emptyset$.
\end{prop}

From now on, we further assume that $G/H \hookrightarrow X_0$ is a spherical embedding
with associated colored fan $\Sigma_0$ such that
$X_0$ only contains $G$-orbits of codimensions $0$ and $1$, \ie we have
\begin{align*}
\Sigma_0 = \{(0, \emptyset)\} \cup \{(\cone(\rho(D)), \emptyset) : D \in \gdiv\}\text{.}
\end{align*}
We continue to denote by $X$ a (normal) variety with the same (graded) Cox ring $R$ as $X_0$.
In particular, if the $G$-action extends to $X$, then 
the variety $X$ is spherical and 
$X_0$ is obtained from $X$ by removing the $G$-orbits of codimension at least $2$.

The description of the divisor class group of a spherical variety
from \cite[Proposition~4.1.1]{brcox} shows that
in the exact sequences
\begin{align*}
\xymatrix@R=0pt{
0 \ar[r] & L_\Q \ar[r] & {(\Q^{\adiv})^*} \ar[r]^-{\Pm} & \Nm_\Q \ar[r] & 0\\
0 & K_\Q \ar[l] & {\phantom{(}\Q^{\adiv}\phantom{)^*}} \ar[l]_-{\Qm} & \Mm_\Q \ar[l] & 0 \ar[l]
}
\end{align*}
of Section~\ref{sec:gdse} we have $K_\Q = \Cl(X_0)_\Q$ and $\Qm(D) = [D]$.
Moreover, Brion has shown that spherical varieties have
finitely generated Cox rings. The following Remark~\ref{rem:sphcox} summarizes the
properties of the Cox ring which we are going to use.

\begin{remark}
\label{rem:sphcox}
It follows from \cite[Theorem~4.3.2]{brcox}, \cite[Theorem~3.6]{g1}, or \cite[Theorem~4.5.4.6]{coxrings} and
\cite[Proposition~2.4]{g1} that there exist
positive integers $n_D$ and elements $f_{D, 1}, \dots, f_{D, n_D} \in R$ with the following
properties:
\begin{enumerate}[(i)]
\item For every $D \in \adiv$ and every $1 \le \ell \le n_D$ we have $\deg f_{D, \ell} = [D]$.
\item The system $\Ff \coloneqq \{f_{D,\ell} : D \in \adiv, 1 \le \ell \le n_D\}$ consists
of pairwise nonassociated $K$-prime generators for $R$ such that the $K$-grading 
is almost free.
\item We have $n_D = 1$ for $D \in \gdiv$ and $n_D \ge 2$ for $D \in \bdiv$.
\item Assume that the $G$-action on $X_0$ extends to $X$.
Then, for every $D \in \adiv$ the (possibly empty) closed subset
$\pi(\V(f_{D,1},\dots,f_{D,n_D})) \subseteq X$ is $G$-stable.
Moreover, for every $G$-orbit $Y \subseteq X$ 
we have $Y \subseteq \overline{D}$ if and only if
$Y \subseteq \pi(\V(f_{D,1},\dots,f_{D,n_D}))$.
\end{enumerate}
\end{remark}

According to \cite[Corollary~3.1.4.6]{coxrings},
the $G$-action on $X_0$ extends to any $X$
arising from a bunched ring $(R, \Ff, \Theta)$
with a maximal bunch $\Theta$. Then, the following
Remark~\ref{rem:str} together with Remark~\ref{rem:sphcox}(iv) shows
that the $G$-action on $X_0$ also extends to $X$ when the bunch
is not maximal.

\begin{remark}
\label{rem:str}
Let $X$ arise from the bunched ring $(R, \Ff, \Theta)$.
Then we have 
\begin{align*}
X(\tau) &= \bigcup_{\substack{J \subseteq \Ff\\\cone(\deg f_{D,\ell} : f_{D,\ell} \in J) = \tau}}\rleft(\bigcap_{f_{D,\ell} \notin J} \pi(\V(f_{D,\ell})) \mathbin{\Big\backslash} \bigcup_{f_{D,\ell}\in J} \pi(\V(f_{D,\ell}))\rright)\\
&= \bigcup_{\substack{J \subseteq \adiv\\\cone(\Qm(J)) = \tau}}
\rleft(\bigcap_{D\notin J}\pi(\V(f_{D,1}, \dots, f_{D,{n_D}})) \mathbin{\Big\backslash} \bigcup_{D\in J} \pi(\V(f_{D,1}, \dots, f_{D,{n_D}})) \rright)
\text{,}
\end{align*}
where the last equality follows from the fact that for every $1 \le \ell \le n_D$ we have $\deg f_{D, \ell} = \Qm(D)$.
\end{remark}

We now explain the relation between the $\Ff$-bunches
and the $\Qm$-bunches from Section~\ref{sec:gdse}.
Note that every projected $\Ff$-face is a $\Qm$-cone.

\begin{lemma}
\label{le:fanbunch}
Let $X$ arise from the bunched ring $(R, \Ff, \Theta)$ and denote by $\Sigma$ the colored fan associated to the spherical embedding $G/H \hookrightarrow X$.
Then for any $I \subseteq \adiv$ such that $(\cone(\Pm(I)), I \cap \bdiv)$ is supported we have $(\cone(\Pm(I)), I \cap \bdiv) \in \Sigma$ if and only if $\cone(\Qm(\adiv \setminus I)) \in \Theta$.
\end{lemma}

\begin{proof}
Let $(\Cm, \Fm) \coloneqq (\cone(\Pm(I)), I \cap \bdiv) \in \Sigma$
and $\tau \coloneqq \cone(\Qm(\adiv \setminus I))$.
Then 
the $G$-orbit corresponding to $(\Cm, \Fm)$ is contained in
$X(\tau)$ by Remarks~\ref{rem:sphcox}(iv) and \ref{rem:str}.
It follows from $X(\tau) \ne \emptyset$ that we have
$\tau \in \Theta$.

Let $\tau \coloneqq \cone(\Qm(\adiv \setminus I)) \in \Theta$ and assume 
that $(\Cm, \Fm) \coloneqq (\cone(\Pm(I)), I \cap \bdiv)$ is supported.
According to Proposition~\ref{prop:qp},
there exists a quasi-projective spherical embedding $G/H \hookrightarrow X'$
with associated colored fan $\Sigma'$ such that $\Sigma_0 \subseteq \Sigma'$,
$(\Cm, \Fm) \in \Sigma'$, and $X' \setminus X_0$ is of codimension at least $2$.
According to Proposition~\ref{prop:ea2}, 
then $X'$ is an open subvariety of a variety arising from a bunched ring $(R, \Ff, \Theta')$, which we call again $X'$.
The $G$-orbit corresponding to $(\Cm, \Fm)$ is then contained
in $X'(\tau)$ by the first part of the proof.
As $X'(\tau) = X(\tau)$, we obtain $(\Cm, \Fm) \in \Sigma$.
\end{proof}

\begin{lemma}
\label{le:suppproj}
The projected $\Ff$-faces are exactly the supported $\Qm$-cones.
\end{lemma}
\begin{proof}
It follows from Lemma~\ref{le:fanbunch} that every projected $\Ff$-face
is a supported $\Qm$-cone. On the other hand, if $\tau$ is a supported $\Qm$-cone,
there exists $(\Cm, \Fm) \in \{\tau\}^\sharp$.
As in the second part of the proof of Lemma~\ref{le:fanbunch},
the $G$-orbit corresponding to $(\Cm, \Fm)$ is contained in some
spherical embedding arising from a bunched ring $(R, \Ff, \Theta)$
with $\tau \in \Theta$. Hence $\tau$ is a projected $\Ff$-face.
\end{proof}

\begin{lemma}
\label{le:fqtrue}
The definitions of \enquote{true} for $\Ff$-bunches
and $\Qm$-bunches coincide.
\end{lemma}
\begin{proof}
For every $D \in \bdiv$ we have $n_D \ge 2$
and $\deg f_{D, \ell} = \Qm(D)$ for every $1 \le \ell \le n_D$.
It follows that we have
\begin{align*}
\cone(\Qm(\adiv)) = \cone(\deg f_{D, \ell} : f_{D, \ell} \in \Ff \setminus \{f_{D_0, \ell_0}\})
\end{align*}
for every $D_0 \in \bdiv$ and every $1 \le \ell_0 \le n_{D_0}$, from which the claim follows.
\end{proof}

\begin{theorem}\label{thm:br}
The true $\Ff$-bunches are exactly the true $\Qm$-bunches. Moreover,
for a true $\Ff$-bunch $\Theta$ the
variety $X$ associated to the bunched ring $(R, \Ff, \Theta)$ is
the spherical embedding $G/H \hookrightarrow X$ associated to the colored fan $\Theta^\sharp$.
\end{theorem}
\begin{proof}
This follows immediately from Lemmas~\ref{le:fanbunch}, \ref{le:suppproj}, and
\ref{le:fqtrue}.
\end{proof}

We can now prove Theorem~\ref{th:m} under the condition $\Gamma(X, \Om^*_X) = \C^*$.

\begin{prop}
Let $G/H \hookrightarrow X$ be a spherical
embedding with $\Gamma(X, \Om^*_X) = \C^*$ and associated colored fan $\Sigma$.
Then $X$ has the $A_2$-property if and
only if any two colored cones in $\Sigma$
intersect in a common face.
\end{prop}
\begin{proof}
According to Proposition~\ref{prop:ea2}, if $X$ has the $A_2$-property,
it is an open subvariety of a variety arising from a bunched ring,
hence any two colored cones in $\Sigma$
intersect in a common face by Theorems~\ref{thm:br} and \ref{thm:fb}

On the other hand, if any two colored cones in $\Sigma$
intersect in a common face, then $\Sigma$ is a true $\Pm$-fan
and can be extended to a true maximal $\Pm$-fan.
It then follows from Theorems~\ref{thm:fb} and \ref{thm:br} that $X$
is an open subvariety of a variety arising from a bunched ring.
Therefore $X$ has the $A_2$-property.
\end{proof}

Under the condition $\Gamma(X, \Om^*_X) = \C^*$,
this proves Theorem~\ref{th:ak} in the case $k=2$
since the $\sharp$-operation sends colored cones in $\Sigma$ which do not intersect
in a common face to cones in $\Sigma^\sharp$ whose relative interiors do not intersect,
while the case $k\ge 3$ follows from Proposition~\ref{prop:eak3}.
The following result shows that the condition $\Gamma(X, \Om^*_X) = \C^*$ can be removed
and completes the proof of Theorems~\ref{th:ak} and \ref{th:m}.

\begin{prop}
Let $G/H \hookrightarrow X$ be a spherical
embedding with associated colored fan $\Sigma$.
Then there exist linearly independent rays $\rho_1, \dots, \rho_d \subseteq \Vm$
such that
\begin{align*}
\Nm_\Q = \vspan_\Q{}\{\rho_1, \dots, \rho_d\} \oplus \vspan_\Q{}\{\Cm : (\Cm, \Fm) \in \Sigma\}\text{.}
\end{align*}
Moreover, for the spherical embedding $G/H \hookrightarrow X'$ associated
to the colored fan $\Sigma' \coloneqq \Sigma \cup \{(\rho_i, \emptyset): 1 \le i \le d\}$,
where we denote by $D_i$ the $G$-invariant prime divisor in $X'$ corresponding to the colored cone $(\rho_i, \emptyset)$,
 the following statements hold:
\begin{enumerate}[(i)]
\item We have $\Gamma(X', \Om_{X'}^*) = \C^*$.
\item We have $\Cl(X')_\Q = \Cl(X)_\Q$ with $[D_i] = 0 \in \Cl(X')_\Q$ for every $1 \le i \le d$\text{.}
\item We have $(\Sigma')^\sharp = \Sigma^\sharp$.
\item $X'$ has the $A_k$-property if and only if $X$ has the $A_k$-property.
\end{enumerate} 
\end{prop}
\begin{proof}
The existence of the rays $\rho_1, \dots, \rho_d$
follows from the fact that the valuation cone $\Vm$ is of
full dimension in $\Nm_\Q$. Then, (i) follows from Remark~\ref{rem:constunits}
and (ii) follows from \cite[Proposition~4.1.1]{brcox}.
From $\Sigma \subseteq \Sigma'$, we obtain $\Sigma^\sharp \subseteq (\Sigma')^\sharp$.
On the other hand, it follows from (ii) that we have
$\{(\rho_i, \emptyset)\}^\sharp = \{(0, \emptyset)\}^\sharp \subseteq \Sigma^\sharp$,
so that we also obtain $(\Sigma')^\sharp \subseteq \Sigma^\sharp$. This proves (iii).

If $X'$ has the $A_k$-property, then the open subvariety $X \subseteq X'$ also
has the $A_k$-property. Now assume that $X'$ does not have the $A_k$-property.
Then there exist $x_1,\dots,x_k \in X'$ which are not contained
in any common affine open neighbourhood.
Since we have proven Theorem~\ref{th:ak} in the case $\Gamma(X', \Om_{X'}^*) = \C^*$,
for every $1 \le i \le k$ we have $x_i \in X'(\tau_i)$ for some $\tau_i \in (\Sigma')^\sharp$
such that $\tau_1^\circ \cap \dots \cap \tau_k^\circ = \emptyset$.
Now let $\tau_0$ be the unique element in $\{(0, \emptyset)\}^\sharp$.
From $\{(\rho_i, \emptyset)\}^\sharp = \{(0, \emptyset)\}^\sharp$,
we obtain $X' \setminus X \subseteq X'(\tau_0)$.
Since $(0, \emptyset)$ is a face of every colored cone in $\Sigma'$,
we have $\tau^\circ \subseteq \tau_0^\circ$ for every
$\tau \in (\Sigma')^\sharp$.
By induction, we may assume that $X'$ does have the $A_{k-1}$-property,
hence it follows from $\tau_1^\circ \cap \dots \cap \tau_k^\circ = \emptyset$
that we have $\tau_i \ne \tau_0$ for every $1 \le i \le k$.
Therefore we have $x_1, \dots, x_k \in X$, so that $X$ does not have the $A_k$-property.
\end{proof}

Using Theorem~\ref{thm:br}, it is possible to apply results on bunched rings
to spherical varieties. We give two examples.

\begin{remark}[The criterion for $\Q$-factoriality]
We recover Remark~2.5 by combining \cite[Corollary~3.3.1.9]{coxrings},
Theorem~\ref{thm:br}, and the last statement of Theorem~\ref{thm:fb}.
\end{remark}

\begin{remark}[{The canonical toric embedding; see \cite[3.2.5]{coxrings}}]
Let $X$ arise from the bunched ring $(R, \Ff, \Theta)$. 
Consider the $K$-graded polynomial ring $\C[\Ff]$
where the elements $f_{D,\ell} \in \Ff$ are interpreted as homogeneous variables
of degree $[D] \in K$.
Let $\Q^{\Ff}$ and $(\Q^{\Ff})^*$ be dual vector
spaces with respective standard bases $\smash{\{e_{f_{D,\ell}} : f_{D,\ell} \in \Ff\}}$ and
$\smash{\{e_{f_{D,\ell}}^* : f_{D,\ell} \in \Ff\}}$, which are dual to each other.
The map $Q\colon \Q^{\Ff} \to K_\Q$ with $e_{f_{D,\ell}} \mapsto \deg f_{D,\ell} = [D]$
induces the following pair of mutually dual exact sequences of vector spaces.
\begin{align*}
\xymatrix@R=0pt{
0 \ar[r] & L_\Q \ar[r] & {(\Q^{\Ff})^*} \ar[r]^-{P} & N_\Q \ar[r] & 0\\
0 & K_\Q \ar[l] & {\phantom{(}\Q^{\Ff}\phantom{)^*}} \ar[l]_-{Q} & M_\Q \ar[l] & 0 \ar[l]
}
\end{align*}
Let $Z$ be the (toric) variety arising from the bunched ring $(\C[\Ff], \Ff, \Theta)$.
Its fan can be obtained as $\Theta^\sharp$, where now $\sharp$ denotes the
$\sharp$-operation with respect to the exact sequences given here, \ie yielding a fan in $N_\Q$.
The surjective homomorphism of graded (Cox) rings $\C[\Ff] \to R$ induces a closed embedding $X \hookrightarrow Z$.
\end{remark}

\section*{Acknowledgements}
The author would like to thank Victor Batyrev for encouragement and advice
as well as J\"{u}rgen Hausen for several highly useful discussions.

\bibliographystyle{amsalpha}
\bibliography{sak}

\end{document}